\documentclass{amsart}
\usepackage{amssymb}

\title[Hall-Littlewood Polynomials of Type $B$ and $C$]{Haglund-Haiman-Loehr Type Formulas for Hall-Littlewood Polynomials of Type $B$ and $C$}

\author{Cristian Lenart}

\address{Department of Mathematics and Statistics, State University of New York at Albany, Albany, NY 12222}
\email{lenart@albany.edu}



\thanks{Cristian Lenart was partially supported by the National Science Foundation grant  DMS-0701044}

\DeclareMathOperator{\des}{des}

\DeclareMathOperator{\inv}{inv}
\DeclareMathOperator{\cinv}{cinv}

\DeclareMathOperator{\rev}{rev}

\newlength{\cellsize}
\newcommand\tableau[1]{
\vcenter{
\let\\=\cr
\baselineskip=-16000pt
\lineskiplimit=16000pt
\lineskip=0pt
\halign{&\tableaucell{##}\cr#1\crcr}}}

\newcommand{\tableaucell}[1]{{%
\def \arg{#1}\def \void{}%
\ifx \void \arg
\vbox to \cellsize{\vfil \hrule width \cellsize height 0pt}%
\else
\unitlength=\cellsize
\begin{picture}(1,1)
\put(0,0){\makebox(1,1){$#1$}}
\put(0,0){\line(1,0){1}}
\put(0,1){\line(1,0){1}}
\put(0,0){\line(0,1){1}}
\put(1,0){\line(0,1){1}}
\end{picture}%
\fi}}

\numberwithin{equation}{section}

\theoremstyle{plain}
\newtheorem{theorem}{Theorem}[section]
\newtheorem{proposition}[theorem]{Proposition}
\newtheorem{lemma}[theorem]{Lemma}

\newtheorem{definition}[theorem]{Definition}

\newtheorem{example}[theorem]{Example}

\theoremstyle{remark}
\newtheorem{remark}[theorem]{Remark}
\newtheorem{remarks}[theorem]{Remarks}

\def\R{\mathbb{R}}
\def\Z{\mathbb{Z}}

\def\F{{\mathcal F}}

\def\Waff{W_{\mathrm{aff}}}

\def\h{\mathfrak{h}}
\def\hR{\mathfrak{h}^*_\mathbb{R}}
\newcommand{\casethree}[6]{\left\{ \begin{array}{ll} #1 &\mbox{if $#2$} \\#3 &\mbox{if $#4$} \\ #5 & \mbox{if $#6$}\,. \end{array} \right.}

\newcommand{\casetwo}[3]{\left\{ \begin{array}{ll} #1 &\mbox{if $#2$} \\ #3 &\mbox{otherwise}\,. \end{array} \right.}

\newcommand{\casetwoex}[4]{\left\{ \begin{array}{ll} #1 &\mbox{if $#2$} \\ #3 &\mbox{if $#4$} \,. \end{array} \right.}
\newcommand{\casetwoexc}[4]{\left\{ \begin{array}{ll} #1 &\mbox{if $#2$} \\ #3 &\mbox{if $#4$} \,, \end{array} \right.}

\newcommand{\stacksum}[2]{\sum_{\begin{array}{c}\vspace{-5.4mm}\;\\ \vspace{-1mm}\scriptstyle{#1}\\ \scriptstyle{#2}\end{array}} }

\begin{document}
\bibliographystyle{plain}

\begin{abstract} In previous work we showed that two apparently unrelated formulas  for the Hall-Littlewood polynomials of type $A$ are, in fact, closely related. The first is the tableau formula obtained by specializing $q=0$ in the Haglund-Haiman-Loehr formula for Macdonald polynomials. The second is the type $A$ instance of Schwer's formula (rephrased and rederived by Ram) for Hall-Littlewood polynomials of arbitrary finite type; Schwer's formula is in terms of so-called alcove walks, which originate in the work of Gaussent-Littelmann and of the author with Postnikov on discrete counterparts to the Littelmann path model. We showed that the tableau formula follows by ``compressing'' Ram's version of Schwer's formula. In this paper, we derive tableau formulas for the Hall-Littlewood polynomials of type $B$ and $C$ by compressing the corresponding instances of Schwer's formula. 
\end{abstract}

\maketitle

\section{Introduction}

Hall-Littlewood polynomials are at the center of many recent developments in representation theory and algebraic combinatorics. They were originally defined in type $A$, as a basis for the algebra of symmetric functions depending on a parameter $t$; this basis interpolates between two fundamental bases: the one of Schur functions, at $t=0$, and the one of monomial functions, at $t=1$. Beside the original motivation for defining Hall-Littlewood polynomials which comes from the Hall algebra \cite{litcsf}, there are many other applications (see e.g. \cite{lenhlp} and the references therein). 

Macdonald \cite{macsfg} showed that there is a formula for the spherical functions corresponding to a Chevalley group over a $p$-adic field which generalizes the formula for the Hall-Littlewood polynomials. Thus, the Macdonald spherical functions generalize the Hall-Littlewood polynomials to all root systems, and the two names are used interchangeably in the literature. There are two families of Hall-Littlewood polynomials of arbitrary type, called $P$-polynomials and $Q$-polynomials, which form dual bases for the Weyl group invariants. The $P$-polynomials specialize to the Weyl characters at $t=0$.  The transition matrix between Weyl characters and $P$-polynomials is given by Lusztig's $t$-analog of weight multiplicities (Kostka-Foulkes polynomials of arbitrary type), which are certain affine Kazhdan-Lusztig polynomials \cite{katsfq,lusscq}. On the combinatorial side, we have the Lascoux-Sch\"utzenberger formula for the Kostka-Foulkes polynomials in type $A$ \cite{lassuc}, but no generalization of this formula to other types is known. Other applications of the type $A$ Hall-Littlewood polynomials that extend to arbitrary type are those related to fermionic multiplicity formulas \cite{aakfpk} and affine crystals \cite{laslqa}. We refer to \cite{narkfp,stekfp} for surveys on Hall-Littlewood polynomials of arbitrary type. 

Macdonald \cite{macsft,macopa} defined a remarkable family of orthogonal polynomials depending on parameters $q,t$, which bear his name. These polynomials generalize the spherical functions for a $p$-adic group, the Jack polynomials, and the zonal polynomials. At $q=0$, the Macdonald polynomials specialize to the Hall-Littlewood polynomials, and thus they further specialize to the Weyl characters (upon setting $t=0$ as well). There has been considerable interest recently in the combinatorics of Macdonald polynomials. This stems in part from a combinatorial formula for the ones corresponding to type $A$, which is due to Haglund, Haiman, and Loehr \cite{hhlcfm}. This formula is in terms of fillings of Young diagrams, and uses two statistics, called inv and maj, on such fillings. The Haglund-Haiman-Loehr formula already found important applications, such as new proofs of the positivity theorem for Macdonald polynomials, which states that the two-parameter Kostka-Foulkes polynomials have nonnegative integer coefficients. One of the mentioned proofs, due to Grojnowski and Haiman \cite{gahaha}, is based on Hecke algebras, while the other, due to Assaf \cite{asasem}, is purely combinatorial and leads to a positive formula for the two-parameter Kostka-Foulkes polynomials. Moreover, in the one-parameter case (i.e., when $q=0$), the Haglund-Haiman-Loehr formula was used to give a concise derivation of the Lascoux-Sch\"utzenberger formula for the Kostka-Foulkes polynomials of type $A$ \cite[Section 7]{hhlcfm}.

An apparently unrelated development, at the level of arbitrary finite root systems, led to Schwer's formula \cite{schghl}, rephrased and rederived by Ram \cite{ramawh}, for the Hall-Littlewood polynomials of arbitrary type. The latter formulas are in terms of so-called alcove walks, which originate in the work of Gaussent-Littelmann \cite{gallsg} and of the author with Postnikov \cite{lapawg,lapcmc} on discrete counterparts to the Littelmann path model \cite{litlrr,litpro}. Schwer's formula was recently generalized by Ram and Yip to a similar formula for the Macdonald polynomials \cite{raycfm}. The generalization consists in the fact that the latter formula is in terms of alcove walks with both ``positive'' and ``negative'' foldings, whereas in the former only ``positive'' foldings appear.

In \cite{lenhlp}, we related Schwer's formula to the Haglund-Haiman-Loehr formula. More precisely, we showed that we can group the terms in the type $A$ instance of Schwer's formula (in fact, we used Ram's version of it) for $P_\lambda(x;t)$ into equivalence classes, such that the sum in each equivalence class is a term in the Haglund-Haiman-Loehr formula for $q=0$. An equivalence class consists of all the terms corresponding to alcove walks that produce the same filling of a Young diagram $\lambda$ (indexing the Hall-Littlewood polynomial) via a simple construction. In fact, we required that the partition $\lambda$ has no two parts identical (i.e., it is a regular weight); the general case, which displays additional complexity, will be considered in a future publication. This work was extended in \cite{lencfm}, by showing that the type $A$ instance of the Ram-Yip formula for Macdonald polynomials compresses, in a similar way, to a formula which is analogous to the Haglund-Haiman-Loehr one, but has fewer terms.

In this paper we extend the results in \cite{lenhlp} to types $B$ and $C$. More precisely, we derive formulas for the Hall-Littlewood polynomials of type $B$ and $C$ indexed by regular weights in terms of fillings of Young diagrams; we do this by compressing the corresponding instances of Schwer's formula (in fact, we again use Ram's version of it). Note that no tableau formula for the Hall-Littlewood or Macdonald polynomials exists beyond type $A$ so far. Our approach provides a natural way to obtain such formulas, and suggests that this method could be further extended to type $D$ (this case is slightly more complex than types $B$ and $C$, see below), as well as to Macdonald polynomials; these problems are currently explored, as is the compression in the case of a Hall-Littlewood polynomial indexed by a non-regular weight. Our formula is more complex than the corresponding one in type $A$ (i.e., the Haglund-Haiman-Loehr formula at $q=0$); however, the statistic we use is, in the case of some special fillings, completely similar to the Haglund-Haiman-Loehr inversion statistic (which is the more intricate of their two statistics). The naturality of our formula is also supported by the fact that the Kashiwara-Nakashima tableaux of type $B$ and $C$ \cite{kancgr} are, essentially, the surviving fillings in this formula when we set $t=0$. We also note that that the passage from (Ram's version of) Schwer's formula to ours results in a considerably larger reduction in the number of terms in type $B$ and $C$ compared to type $A$. 
In terms of applications, it would be very interesting to see whether our formula could be used to derive, in the spirit of \cite[Section 7]{hhlcfm}, a positive combinatorial formula for Lusztig's $t$-analog of weight multiplicities in type $B$ and $C$, which has been long sought.

\section{The tableau formula in type $C$}\label{sectform}

Let us start by recalling the Weyl group of type $B/C$, viewed as the group of signed permutations $B_n$. Such permutations are bijections $w$ from $[\overline{n}]:=\{1<2<\ldots<n<\overline{n}<\overline{n-1}<\ldots<1\}$ to $[\overline{n}]$ satisfying $w(\overline{\imath})=\overline{w(i)}$. Here $\overline{\imath}$ is viewed as $-i$, so $\overline{\overline{\imath}}=i$. We use the window notation $w=w(1)\ldots w(n)$. Given $1\le i<j\le n$, we denote by $(i,j)$ the reflection which transposes the entries in positions $i$ and $j$ (upon right multiplication). Similarly, we denote by $(i,\overline{\jmath})$, again for $i<j$, the transposition of entries in positions $i$ and $j$ followed by the sign change of those entries. Finally, we denote by $(i,\overline{\imath})$ the sign change in position $i$. Given $w$ in $B_n$, we define
\begin{align}\label{defl}&\ell_+(w):=|\{(k,l)\::\:1\le k<l\le n,\,w(k)>w(l)\}|\,,\\&\ell_-(w):=|\{(k,l)\::\:1\le k\le l\le n,\,w(k)>\overline{w(l)}\}|\,.\nonumber\end{align}
Then the length of $w$ is given by $\ell(w):=\ell_+(w)+\ell_-(w)$. 

Let $\lambda$ be a partition corresponding to a regular weight in type $C_n$ for $n\ge 2$, that is $\lambda=(\lambda_1>\lambda_2>\ldots>\lambda_n>0)$. Consider the shape $\widehat{\lambda}$ obtained from $\lambda$ by replacing each column of height $k$ with $k$ or $2k-1$ (adjacent) copies of it, depending on the given column being the first one or not. We are representing a filling $\sigma$ of $\widehat{\lambda}$ as a concatenation of columns $C_{ij}$ and $C'_{ik}$, where $i=1,\ldots,\lambda_1$, while for a given $i$ we have $j=1,\ldots,\lambda_i'$ if $i>1$, $j=1$ if $i=1$, and $k=2,\ldots,\lambda_i'$; the columns $C_{ij}$ and $C_{ik}'$ have height $\lambda_i'$. More precisely, we let
\begin{equation}\label{filling}\sigma={\mathcal C}^{\lambda_1}\ldots {\mathcal C}^{1}\,,\end{equation}
where
\[{\mathcal C}^{i}:=\casetwoex{C_{i2}'\ldots C_{i,\lambda_i'}'C_{i1}\ldots C_{i,\lambda_i'}}{i>1}{C_{i2}'\ldots C_{i,\lambda_i'}'C_{i1}}{i=1}\]
Note that the leftmost column is  $C_{\lambda_1,1}$, and the rightmost column is $C_{11}$. 

Essentially, the above description says that the column to the right of $C_{ij}$ is $C_{i,j+1}$, whereas the column to the right of $C_{ik}'$ is $C_{i,k+1}'$. Here we are assuming that the mentioned columns exist, up to the following conventions:
\begin{equation}\label{conv}C_{i,\lambda_i'+1}=\casetwoexc{C_{i-1,2}'}{i>1\:\mbox{ and }\:\lambda_{i-1}'>1}{C_{i-1,1}}{i>1\:\mbox{ and }\:\lambda_{i-1}'=1}\;\;\;\;\;C_{i,\lambda_i'+1}'=C_{i1}\,.\end{equation}
We consider the set ${\mathcal F}(\lambda)$ of fillings of $\widehat{\lambda}$ with entries in $[\overline{n}]$ which satisfy the following conditions:
\begin{enumerate}
\item the rows are weakly decreasing from left to right;
\item no column contains two entries $a,b$ with $a=\pm b$;
\item each column (with the exception of the leftmost one) is related to its left neighbor as indicated below; essentially, it differs from this neighbor by a ``signed cycle'', that is, a composition $(r_1,\overline{\jmath})\ldots (r_p,\overline{\jmath})$, where $1\le r_1<\ldots <r_p<j$; furthermore, $j$ varies from 1 to the length of the corresponding column, as we consider the columns from left to right.
\end{enumerate}

Here we let reflections in $B_n$ act on columns $C$ like they do on signed permutations; for instance, $C(a,\overline{b})$ is the column obtained from $C$ by transposing the entries in positions $a,b$ and by changing their signs. Let us first explain the passage from some column $C_{ij}$ to $C_{i,j+1}$. There exist positions $1\le r_1<\ldots <r_p<j$ (possibly $p=0$) such that $C_{i,j+1}$ differs from $D=C_{ij}(r_1,\overline{\jmath})\ldots (r_p,\overline{\jmath})$ only in position $j$, while  $C_{i,j+1}(j)\not\in \{\pm D(r)\::\:r\in[\lambda_i']\setminus\{j\}\}$ and $C_{i,j+1}(j)\le D(j)$. To include the case $j=\lambda_i'$ in this description, just replace $C_{i,j+1}$ everywhere by $C_{i,j+1}[1,\lambda_i']$ and use the conventions (\ref{conv}). Let us now explain the passage from some column $C_{ik}'$ to $C_{i,k+1}'$. There exist positions $1\le r_1<\ldots <r_p<k$ (possibly $p=0$) such that $C_{i,k+1}'=C_{ik}'(r_1,\overline{k})\ldots (r_p,\overline{k})$. This description includes the case $k=\lambda_i'$, based on the conventions (\ref{conv}).

Let us now define the content of a filling. For this purpose, we first associate with a filling $\sigma$ a ``compressed'' version of it, namely the filling $\overline{\sigma}$ of the partition $2\lambda$. This is defined as follows:
\begin{equation}\label{redfilling}
\overline{\sigma}=\overline{{\mathcal C}}^{\lambda_1}\ldots \overline{{\mathcal C}}^{1}\,,\;\;\;\;\mbox{where $\:\overline{{\mathcal C}}^{i}:=C_{i2}'C_{i1}$}\,,\end{equation}
where the conventions (\ref{conv}) are used again. Now define ${\rm ct}(\sigma)=(c_1,\ldots,c_n)$, where $c_i$ is half the difference between the number of occurences of the entries $i$ and $\overline{\imath}$ in $\overline{\sigma}$. Sometimes, this vector is written in terms of the coordinate vectors $\varepsilon_i$, as 
\begin{equation}\label{defct}{\rm ct}(\sigma)=c_1\varepsilon_1+\ldots+c_n\varepsilon_n=\frac{1}{2}\sum_{b\in\overline{\sigma}}\varepsilon_{\overline{\sigma}(b)}\,;\end{equation}
here the last sum is over all boxes $b$ of $\overline{\sigma}$, and we set $\varepsilon_{\overline{\imath}}:=-\varepsilon_i$. 

We now define two statistics on fillings that will be used in our compressed formula for Hall-Littlewood polynomials. Intervals refer to the discrete set $[\overline{n}]$. Let 
\begin{equation}\label{defsigma}\sigma_{ab}:=\casetwo{1}{a,b\ge\overline{n}}{0}\end{equation}
Given  a sequence of integers $w$, we write $w[i,j]$ for the subsequence $w(i)w(i+1)\ldots w(j)$. We use the notation $N_{ab}(w)$ for the number of entries $w(i)$ with $a<w(i)<b$. 

Given two columns $D,C$ of the same height $d$ such that $D\ge C$ componentwise, we will define two statistics $N(D,C)$ and $\des(D,C)$ in some special cases, as specified below.

\emph{Case} 0. If $D=C$, then $N(D,C):=0$ and $\des(D,C):=0$. 

\emph{Case} 1. Assume that $C=D(r,\overline{\jmath})$ with $r<j$. Let $a:=D(r)$ and $b:=D(j)$. In this case, we set
\[N(D,C):=N_{\overline{b}a}(D[r+1,j-1])+|(\overline{b},a)\setminus\{\pm D(i)\::\:i=1,\ldots,j\}|+\sigma_{ab}\,,\]
and $\des(D,C):=1\,.$

\emph{Case} 2. Assume that $C=D(r_1,\overline{\jmath})\ldots (r_p,\overline{\jmath})$ where $1\le r_1<\ldots <r_p<j$. Let $D_i:=D(r_1,\overline{\jmath})\ldots (r_i,\overline{\jmath})$ for $i=0,\ldots,p$, so that $D_0=D$ and $D_p=C$. We define
\[N(D,C):=\sum_{i=1}^p N(D_{i-1},D_i)\,,\;\;\;\des(D,C):=p\,.\]

\emph{Case} 3. Assume that $C$ differs from $D':=D(r_1,\overline{\jmath})\ldots (r_p,\overline{\jmath})$ with $1\le r_1<\ldots <r_p<j$ (possibly $p=0$) only in position $j$, while  $C(j)\not\in \{\pm D'(r)\::\:r\in[d]\setminus\{j\}\}$ and $C(j)< D'(j)$. We define
\[N(D,C):=N(D,D')+N_{C(j),D'(j)}(D[j+1,d])\,,\;\;\;\des(D,C):=p+1\,.\]

If the height of $C$ is larger than the height $d$ of $D$ (necessarily by 1), and $N(D,C[1,d])$ can be computed as above, we let $N(D,C):=N(D,C[1,d])$ and $\des(D,C):=\des(D,C[1,d])$. Given a filling $\sigma$ in ${\mathcal F}(\lambda)$ with columns $C_m,\ldots ,C_1$, we set
\[N(\sigma):=\sum_{i=1}^{m-1}N(C_{i+1},C_i)+\ell_+(C_1)\,;\]
here $\ell_+(C_1)$ is defined like in (\ref{defl}). Furthermore, we also set 
 \[\des(\sigma):=\sum_{i=1}^{m-1}\des(C_{i+1},C_i)\,.\]
Note that $\des(\sigma)$ essentially counts the descents in the rows of $\sigma$. 

We can now state our new formula for the Hall-Littlewood polynomials of type $C$, which follows as a corollary of our main result, i.e., Theorem \ref{ccompress}. A completely similar formula in type $B$ is discussed in Section \ref{typeb}. We refer to Proposition \ref{relninv} and Remarks \ref{specialfill} for more insight into our formula. In particular, note that the Kashiwara-Nakashima tableaux of type $C$ are, essentially, the surviving fillings in this formula when we set $t=0$, and that, in some special cases, the statistic $N(\sigma)$ is completely similar to the Haglund-Haiman-Loehr inversion statistic (the more intricate of their two statistics). 

\begin{theorem}\label{newformc}
Given a regular weight $\lambda$, we have
\begin{equation}\label{newform}P_\lambda(X;t)=\sum_{\sigma\in{\mathcal F}(\lambda)} t^{N(\sigma)}\,(1-t)^{\des(\sigma)}\,x^{\rm ct(\sigma)}\,,\end{equation}
where $x^{(c_1,\ldots,c_n)} := x_1^{c_1}\ldots x_n^{c_n}$. 
\end{theorem}

\begin{example}{\rm  Consider the simplest case, namely $n=2$ and $\lambda=(2,1)$. This leads to considering fillings of the shape $(3,2)$ with elements in $[\overline{2}]$, namely
\[\tableau{{e}&{c}&{a}\\&{d}&{b}}\,.\]
The fillings need to satisfy the following conditions:
\begin{enumerate}
\item $a\le c\le e$, $b\le d$;
\item $a\ne\pm b$;
\item either $c=a$ and $d=b$, or $c=\overline{b}$ and $d=\overline{a}$.
\end{enumerate}
For $i\in\{1,2\}$, let $n_i$ be half the difference between the number of $i$'s and $\overline{\imath}$'s in the multiset $\{a,b,c,d,e,e\}$. Given a proposition $A$, we let $\chi(A)$ be 1 or 0, depending on the logical value of $A$ being true or false. Then
\[P_{(2,1)}(x_1,x_2;t)=\sum_{(a,b,c,d,e)} t^{\chi(a>b)+\chi(a,b\le 2,\,a\ne c)}(1-t)^{\chi(a\ne c)+\chi(c\ne e)}x_1^{n_1}x_2^{n_2}\,.\]
It turns out that there are 27 terms in this sum, versus 70 terms in (Ram's version of) Schwer's formula. For instance, the terms contributing to the coefficient of $x_2$ correspond to the fillings
\[\tableau{{\overline{1}}&{1}&{1}\\&{2}&{2}}\,,\;\;\;\;\;\tableau{{\overline{1}}&{2}&{2}\\&{1}&{1}}\,,\;\;\;\;\;\tableau{{2}&{2}&{1}\\&{\overline{1}}&{\overline{2}}}\,;\]
the associated polynomials in $t$ are 
\[1-t\,,\;\;\;\;\; t(1-t)\,,\;\;\;\;\;1-t\,,\]
respectively. Note that these polynomials are obtained by compressing 3, 2, and 2 terms in Schwer's formula, respectively. By symmetry, the coefficients of $x_1$, $x_2$, $x_1^{-1}$, and $x_2^{-1}$ in $P_{(2,1)}(x_1,x_2;t)$ are all $(t+2)(1-t)$. Other fillings have an even larger number of terms in Schwer's formula corresponding to them, such as
\[\tableau{{\overline{1}}&{\overline{2}}&{\overline{2}}\\&{\overline{1}}&{\overline{1}}}\,,\]
which has 7; in other words, the associated polynomial in $t$, namely $1-t$, which contributes to the coefficient of $x_1^{-2}x_2^{-1}$, is the sum of 7 polynomials of the form $t^r(1-t)^s$ in Schwer's formula. In conclusion, we have 
\begin{align*}P_{(2,1)}(x_1,x_2;t)&=x_1^2x_2+x_1x_2^2+x_1^2x_2^{-1}+x_1x_2^{-2}+x_1^{-1}x_2^2+x_1^{-2}x_2+\\&+x_1^{-1}x_2^{-2}+x_1^{-2}x_2^{-1}+(t+2)(1-t)(x_1+x_2+x_1^{-1}+x_2^{-1})\,.\end{align*}
}
\end{example}

In order to relate our statistic $N(\sigma)$ to the Haglund-Haiman-Loehr {\em inversion statistic}, let us recall some definitions from \cite{hhlcfm,lenhlp}. We start by considering fillings $\tau$ of the shape $\lambda$ with entries in $[\overline{n}]$, which are displayed in ``Japanese style'', as a sequence of columns $\tau=C_{\lambda_1}\ldots C_1$; here $C_i=(C_i(1),\ldots,C_i(\lambda_i'))$, so the entry in cell $u=(i,j)$ is $\tau(u)=C_j(i)$. Two cells $u,v\in \lambda $ are said to {\it attack} each other if they are in one of the following two relative positions:
\begin{equation*}
\tableau{{\bullet}&{}\\ {}&{}\\{\bullet}&{}\\&{}} \, ,\qquad
\tableau{{}&{\bullet}\\ {}&{}\\{\bullet}&{}\\&{}}\; .
\end{equation*}    
An {\it
inversion} of $\tau $ is a pair of attacking cells $(u,v)$ which have one of the following two relative positions, where $a:=\tau(u)<b:=\tau(v)$:
\begin{equation*}
\tableau{{a}&{}\\ {}&{}\\{b}&{}\\&{}} \, ,\qquad 
\tableau{{}&{b}\\ {}&{}\\{a}&{}\\&{}}\,.
\end{equation*}
The Haglund-Haiman-Loehr statistic $\inv(\tau)$ is defined as the number of inversions of $\tau$. The {\em descent statistic}, denoted ${\des}(\tau)$ (which is similar to $\des$ for fillings of $\widehat{\lambda}$ defined above, see below), is the number of cells $u=(i,j)$ with $j\ne 1$ and $\tau(u)>\tau(v)$, where $v=(i,j-1)$. Now let, as usual, $n(\lambda):=\sum_i(i-1)\lambda_i$, and assume that $\tau$ has the following two properties: (i) $\tau(u)\ne\tau(v)$ whenever  $u$ and $v$ attack each other; (ii) $\tau$ is weakly decreasing in rows. Then it was shown in \cite{lenhlp}[Proposition 2.12] that the so-called {\em complementary inversion statistic} $\cinv(\tau):=n(\lambda)-\inv(\tau)$ counts the  triples of cells filled with $a<b<c$  which have the following relative position (here the third cell might be outside the shape $\lambda$, in which case we only require $a<b$):
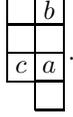
\[\tableau{{}&{b}\\ {}&{}\\{c}&{a}\\&{}}\,.\]

\begin{proposition}\label{relninv} Let $\sigma$ in ${\mathcal F}(\lambda)$ be a filling satisfying the following properties: {\rm (1)} $C_{i,j+1}'=C_{i,j}'$ for all $i$ and $j=1,\ldots,\lambda_i'$; {\rm (2)} $C_{i,j+1}$ only differs from $C_{ij}$ in position $j$. Let $\widetilde{\sigma}$ be the filling of $\lambda$ given by
\[\widetilde{\sigma}:=C_{\lambda_1,1}C_{\lambda_1-1,1}\ldots C_{11}\,.\]
Then $N(\sigma)=\cinv(\widetilde{\sigma})$ and $\des(\sigma)={\des}(\widetilde{\sigma})$. 
\end{proposition}

\begin{proof}
The equality $\des(\sigma)={\des}(\widetilde{\sigma})$ is clear, so we concentrate on the first one. Let $m:=\lambda_1$ be the number of columns of $\lambda$, and let $C_m=C_{m1},\ldots,C_1=C_{11}$ be the columns of $\widetilde{\sigma}$, of lengths $c_m:=\lambda_m',\ldots,c_1:=\lambda_1'$; let $C_k':=C_k[1,c_{k+1}]$, for $k=1,\ldots,m-1$. We refer to a pair $(i,j)$ with $1\le i<j\le c_k$ and $C_k(i)>C_k(j)$ as a (type $A$) inversion in $C_k$. It is easy to see that $\widetilde{\sigma}$ satisfies the properties considered above: (i) $\widetilde{\sigma}(u)\ne\widetilde{\sigma}(v)$ whenever  $u$ and $v$ attack each other; (ii) $\widetilde{\sigma}$ is weakly decreasing in rows. 
We start by evaluating $N({\mathcal C}^kC_{k-1,1})$ (see (\ref{filling})). By definition, we have
\[N({\mathcal C}^kC_{k-1,1})=\sum_{i=1}^{c_k-1} N_{C_{k-1}(i),C_k(i)}(C_k[i+1,c_k])\,.\]
This is the number of inversions $(i,j)$ in $C_k$  for which $C_{k-1}(i)<C_k(j)$. If $(i,j)$ is an inversion in $C_k$ not satisfying the previous condition, then $C_{k-1}(i)> C_k(j)$ (by property (i) of $\widetilde{\sigma}$), and thus $(i,j)$ is an inversion in $C_{k-1}'$ (by property (ii) of $\widetilde{\sigma}$). Moreover, the only inversions of $C_{k-1}'$ which do not arise in this way are those counted by the statistic $\cinv(C_kC_{k-1}')$, so 
\[N({\mathcal C}^kC_{k-1,1})=\ell_+(C_k)-(\ell_+(C_{k-1}')-\cinv(C_kC_{k-1}'))\,.\]

We conclude that
\[N(\sigma)-\ell_+(C_1)=\sum_{k=2}^m \ell_+(C_k)-\ell_+(C_{k-1}')+\cinv(C_kC_{k-1}')\,.\]
Now recall that $\lambda$ has no two parts identical. We clearly have $c_m=1$ so $\ell_+(C_m)=0$. Therefore, we have
\begin{align*}N(\sigma)&=\sum_{k=2}^{m} \ell_+(C_{k-1})-\ell_+(C_{k-1}')+\cinv(C_kC_{k-1}')=\\
&=\sum_{k=2}^{m} \cinv(C_kC_{k-1})=\cinv(\widetilde{\sigma})\,.\end{align*}
\end{proof}


\section{Background on Ram's version of Schwer's formula}\label{prelim}

We recall some background information on finite root systems and affine Weyl groups.

\subsection{Root systems}\label{rootsyst}

Let $\mathfrak{g}$ be a complex semisimple Lie algebra, and $\h$ a Cartan subalgebra, whose rank is $r$.
Let $\Phi\subset \h^*$ be the 
corresponding irreducible {\it root system}, $\hR\subset \h^*$ the real span of the roots, and $\Phi^+\subset \Phi$ the set of positive roots. 
Let $\alpha_1,\ldots,\alpha_r\in\Phi^+$ be the corresponding 
{\it simple roots}.
We denote by $\langle\,\cdot\,,\,\cdot\,\rangle$ the nondegenerate scalar product on $\hR$ induced by
the Killing form.  
Given a root $\alpha$, we consider the corresponding {\it coroot\/} $\alpha^\vee := 2\alpha/\langle\alpha,\alpha\rangle$ and reflection $s_\alpha$.  

Let $W$ be the corresponding  {\it Weyl group\/}, whose Coxeter generators are denoted, as usual, by $s_i:=s_{\alpha_i}$. The length function on $W$ is denoted by $\ell(\,\cdot\,)$. The {\em Bruhat graph} on $W$ is the directed graph with edges $u\rightarrow w$ where $w=u s_{\beta}$ for some $\beta\in\Phi^+$, and $\ell(w)>\ell(u)$; we usually label such an edge by $\beta$ and write $u\stackrel{\beta}\longrightarrow w$. The {\em reverse Bruhat graph} is obtained by reversing the directed edges above. The {\em Bruhat order} on $W$ is the transitive closure of the relation corresponding to the Bruhat graph.

The {\it weight lattice\/} $\Lambda$ is given by
\begin{equation}
\Lambda:=\{\lambda\in \hR \::\: \langle\lambda,\alpha^\vee\rangle\in\Z
\textrm{ for any } \alpha\in\Phi\}.
\label{eq:weight-lattice}
\end{equation}
The weight lattice $\Lambda$ is generated by the 
{\it fundamental weights\/}
$\omega_1,\ldots,\omega_r$, which form the dual basis to the 
basis of simple coroots, i.e., $\langle\omega_i,\alpha_j^\vee\rangle=\delta_{ij}$.
The set $\Lambda^+$ of {\it dominant weights\/} is given by
$$
\Lambda^+:=\{\lambda\in\Lambda \::\: \langle\lambda,\alpha^\vee\rangle\geq 0
\textrm{ for any } \alpha\in\Phi^+\}.
$$
The subgroup of $W$ stabilizing a weight $\lambda$ is denoted by $W_\lambda$, and the set of minimum coset representatives in $W/W_\lambda$ by $W^\lambda$. Let $\Z[\Lambda]$ be the group algebra of the weight lattice $\Lambda$, which  has
a $\Z$-basis of formal exponents $\{x^\lambda \::\: \lambda\in\Lambda\}$ with
multiplication $x^\lambda\cdot x^\mu := x^{\lambda+\mu}$.

Given  $\alpha\in\Phi$ and $k\in\Z$, we denote by $s_{\alpha,k}$ the reflection in the affine hyperplane
\begin{equation}
H_{\alpha,k} := \{\lambda\in \hR \::\: \langle\lambda,\alpha^\vee\rangle=k\}.
\label{eqhyp}
\end{equation}
These reflections generate the {\it affine Weyl group\/} $\Waff$ for the {\em dual root system} 
$\Phi^\vee:=\{\alpha^\vee \::\: \alpha\in\Phi\}$. 
The hyperplanes $H_{\alpha,k}$ divide the real vector space $\hR$ into open
regions, called {\it alcoves.} 
The {\it fundamental alcove\/} $A_\circ$ is given by 
$$
A_\circ :=\{\lambda\in \hR \::\: 0<\langle\lambda,\alpha^\vee\rangle<1 \textrm{ for all }
\alpha\in\Phi^+\}.
$$

\subsection{Alcove walks}\label{alcovewalks}

We say that two alcoves $A$ and $B$ are {\it adjacent} 
if they are distinct and have a common wall.  
Given a pair of adjacent alcoves $A\ne B$ (i.e., having a common wall), we write 
$A\stackrel{\beta}\longrightarrow B$ if the common wall 
is of the form $H_{\beta,k}$ and the root $\beta\in\Phi$ points 
in the direction from $A$ to $B$.  

\begin{definition}
An {\em alcove path\/} is a sequence of alcoves
 such that any two consecutive ones are adjacent. 
We say that an alcove path $(A_0,A_1,\ldots,A_m)$ is {\it reduced\/} if $m$ is the minimal 
length of all alcove paths from $A_0$ to $A_m$.
\end{definition}

We need the following generalization of alcove paths.

\begin{definition}\label{defalcwalk} An {\em alcove walk} is a sequence 
$\Omega=(A_0,F_1,A_1, F_2, \ldots , F_m, A_m, F_{\infty})$ 
such that $A_0,\ldots,$ $A_m$ are alcoves; 
$F_i$ is a codimension one common face of the alcoves $A_{i-1}$ and $A_i$,
for $i=1,\ldots,m$; and 
$F_{\infty}$ is a vertex of the last alcove $A_m$. The weight $F_\infty$ is called the {\em weight} of the alcove walk, and is denoted by $\mu(\Omega)$. 
\end{definition}
 
The {\em folding operator} $\phi_i$ is the operator which acts on an alcove walk by leaving its initial segment from $A_0$ to $A_{i-1}$ intact and by reflecting the remaining tail in the affine hyperplane containing the face $F_i$. In other words, we define
$$\phi_i(\Omega):=(A_0, F_1, A_1, \ldots, A_{i-1}, F_i'=F_i,  A_{i}', F_{i+1}', A_{i+1}', \ldots,  A_m', F_{\infty}')\,;$$
here $A_j' := \rho_i(A_j)$ for $j\in\{i,\ldots,m\}$, $F_j':=\rho_i(F_j)$ for $j\in\{i,\ldots,m\}\cup\{\infty\}$, and $\rho_i$ is the affine reflection in the hyperplane containing $F_i$. Note that any two folding operators commute. An index $j$ such that $A_{j-1}=A_j$ is called a {\em folding position} of $\Omega$. Let $\mbox{fp}(\Omega):=\{ j_1<\ldots< j_s\}$ be the set of folding positions of $\Omega$. If this set is empty, $\Omega$ is called {\em unfolded}. Given this data, we define the operator ``unfold'', producing an unfolded alcove walk, by
\[\mbox{unfold}(\Omega)=\phi_{j_1}\ldots \phi_{j_s} (\Omega)\,.\]

\begin{definition} An  alcove walk
$\Omega=(A_0,F_1,A_1, F_2, \ldots , F_m, A_m, F_{\infty})$ is called {\em positively folded} if, for any folding position $j$, the alcove $A_{j-1}=A_j$ lies on the positive side of the affine hyperplane containing the face $F_j$.
\end{definition} 

We now fix a dominant weight $\lambda$ and a reduced alcove path $\Pi:=(A_0,A_1,\ldots,A_m)$ from $A_\circ=A_0$ to its translate $A_\circ + \lambda=A_m$. Assume that we have
\[A_0\stackrel{\beta_1}\longrightarrow A_1\stackrel{\beta_2}\longrightarrow \ldots
\stackrel{\beta_m}\longrightarrow A_{m}\,,\]
where $\Gamma:=(\beta_1,\ldots,\beta_m)$ is a sequence of positive roots. This sequence, which determines the alcove path, is called a {\em $\lambda$-chain} (of roots). Two equivalent definitions of $\lambda$-chains (in terms of reduced words in affine Weyl groups, and an interlacing condition) can be found in \cite{lapawg}[Definition 5.4] and \cite{lapcmc}[Definition 4.1 and Proposition 4.4]; note that the $\lambda$-chains considered in the mentioned papers are obtained by reversing the ones in the present paper. We also let $r_i:=s_{\beta_i}$, and let $\widehat{r}_i$ be the affine reflection in the common wall of $A_{i-1}$ and $A_i$, for $i=1,\ldots,m$; in other words, $\widehat{r}_i:=s_{\beta_i,l_i}$, where $l_i:=|\{j\le i\::\: \beta_j = \beta_i\}|$ is the cardinality of the corresponding set. Given $J=\{j_1<\ldots<j_s\}\subseteq[m]:=\{1,\ldots,m\}$, we define the Weyl group element $\phi(J)$ and the weight $\mu(J)$ by
\begin{equation}\label{defphimu}\phi(J):={r}_{j_1}\ldots {r}_{j_s}\,,\;\;\;\;\;\mu(J):=\widehat{r}_{j_1}\ldots \widehat{r}_{j_s}(\lambda)\,.\end{equation}

 Given $w\in W$, we define the alcove path $w(\Pi):=(w(A_0),w(A_1),\ldots,w(A_m))$. Consider the set of alcove paths
\[{\mathcal P}(\Gamma):=\{w(\Pi)\::\:w\in W^\lambda\}\,.\]
We identify any $w(\Pi)$ with the obvious unfolded alcove walk of weight $\mu(w(\Pi)):=w(\lambda)$. Let us now consider the set of alcove walks
\[{\mathcal F}_+(\Gamma):=\{\,\mbox{positively folded alcove walks $\Omega$}\::\:\mbox{unfold}(\Omega)\in{\mathcal P}(\Gamma)\}\,.\]
We can encode an alcove walk $\Omega$ in ${\mathcal F}_+(\Gamma)$ by the pair $(w,J)$ in $W^\lambda\times 2^{[m]}$, where 
\[\mbox{fp}(\Omega)=J\;\;\;\;\mbox{and}\;\;\;\;\mbox{unfold}(\Omega)=w(\Pi)\,.\]
Clearly, we can recover $\Omega$ from $(w,J)$ with $J=\{j_1<\ldots<j_s\}$ by
\[\Omega=\phi_{j_1}\ldots \phi_{j_s} (w(\Pi))\,.\]
Let ${\mathcal A}(\Gamma)$ be the image of ${\mathcal F}_+(\Gamma)$ under the map $\Omega\mapsto (w,J)$. We call a pair $(w,J)$ in ${\mathcal A}(\Gamma)$ an {\em admissible pair}, and the subset $J\subseteq[m]$ in this pair a $w$-{\em admissible subset}. The following straightforward result is taken from \cite{lenhlp}.

\begin{proposition}\cite{lenhlp}\label{admpairs} {\rm (1)} We have
\begin{align}\label{decch}{\mathcal A}(\Gamma)=\{\,&(w,J)\in W^\lambda\times 2^{[m]}\::\: J=\{j_1<\ldots<j_s\}\,,\\&w>wr_{j_1}>\ldots>wr_{j_1}\ldots r_{j_s}=w\phi(J)\}\,;\nonumber\end{align}
here the decreasing chain is in the Bruhat order on the Weyl group, its steps not being covers necessarily.

{\rm (2)} If $\Omega\mapsto (w,J)$, then
\[\mu(\Omega)=w(\mu(J))\,.\] 
\end{proposition}

The formula for the Hall-Littlewood $P$-polynomials in \cite{schghl} was rederived in \cite{ramawh} in a slightly different version, based on positively folded alcove walks. Based on Proposition \ref{admpairs}, we now restate the latter formula in terms of admissible pairs. 

\begin{theorem}\cite{ramawh,schghl} \label{hlpthm} Given a dominant weight $\lambda$, we have
\begin{equation}\label{hlpform}P_{\lambda}(X;t)=\sum_{(w,J)\in{\mathcal A}(\Gamma)}t^{\frac{1}{2}(\ell(w)+\ell(w\phi(J))-|J|)}\,(1-t)^{|J|}\,x^{w(\mu(J))}\,.\end{equation}
\end{theorem}

\section{Specializing Ram's version of Schwer's formula to type $C$}\label{specschwer} We now restrict ourselves to the root system of type $C_n$. 
We can identify the space $\h_\R^*$ with  
$V:=\R^n$, the coordinate vectors being $\varepsilon_1,\ldots,\varepsilon_n$.  
The root system $\Phi$ can be represented as 
$\Phi=\{\pm\varepsilon_i\pm\varepsilon_j \::\:  1\leq i<j\leq n\}\cup\{\pm 2\varepsilon_i\::\: 1\leq i\leq n\}$. 
The simple roots are $\alpha_i=\varepsilon_i-\varepsilon_{i+1}$, 
for $i=1,\ldots,n-1$ and $\alpha_n=2\varepsilon_n$. 
The fundamental weights are $\omega_i = \varepsilon_1+\ldots +\varepsilon_i$, 
for $i=1,\ldots,n$. 
The weight lattice is $\Lambda=\Z^n$. A dominant weight $\lambda=\lambda_1\varepsilon_1+\ldots+\lambda_{n-1}\varepsilon_{n-1}+\lambda_n\varepsilon_n$ is identified with the partition $(\lambda _{1}\geq \lambda _{2}\geq \ldots \geq \lambda _{n-1}\geq\lambda_n\geq 0)$ of length at most $n$. A dominant weight is regular if all the previous inequalities are strict (i.e., the corresponding partition has all parts distinct and nonzero). We fix such a partition $\lambda$ for the remainder of this paper.

The corresponding Weyl group $W$ is the group of signed permutations $B_n$. For simplicity, we use the same notation for roots and the corresponding reflections, cf. Section \ref{sectform}. For instance, given $1\le i<j\le n$, we denote by $(i,j)$ the positive root $\varepsilon_i-\varepsilon_j$, by $(i,\overline{\jmath})$ the positive root $\varepsilon_i+\varepsilon_j$, and by $(i,\overline{\imath})$ the positive root $2\varepsilon_i$. 

Let 
\[\Gamma(k):=\Gamma_2'\ldots\Gamma_k'\Gamma_1(k)\ldots\Gamma_k(k)\,,\]
where
\begin{align*}
&\;\:\Gamma_j':=((1,\overline{\jmath}),(2,\overline{\jmath}),\ldots,(j-1,\overline{\jmath}))\,,\\
&\begin{array}{llllll}\Gamma_j(k):=(\!\!\!\!\!&(1,\overline{\jmath}),&(2,\overline{\jmath}),&\ldots,&(j-1,\overline{\jmath}),\\
&(j,\overline{k+1}),&(j,\overline{k+2}),&\ldots,&(j,\overline{n}),\\&(j,\overline{\jmath}),\\ 
&(j,n),&(j,n-1),&\ldots,&(j,k+1)\,)\,.\end{array}\end{align*}

\begin{lemma} $\Gamma(k)$ is an $\omega_k$-chain.
\end{lemma}

\begin{proof} We use the criterion for $\lambda$-chains in \cite{lapcmc}[Definition 4.1, Proposition 4.4], cf. \cite{lapcmc}[Proposition 10.2]. This criterion says that a chain of roots $\Gamma$ is a $\lambda$-chain if and only if it satisfies the following conditions:
\begin{enumerate}
\item[(R1)] The number of occurrences of any positive root $\alpha$ in $\Gamma$ is 
$\langle \lambda,\alpha^\vee \rangle$. 
\item[(R2)] For each triple of positive roots $(\alpha,\,\beta,\,\gamma)$ with $\gamma^\vee=\alpha^\vee+\beta^\vee$, the subsequence of $\Gamma$ consisting of $\alpha,\,\beta,\,\gamma$ is a concatenation of pairs $(\gamma,\alpha)$ and $(\gamma,\beta)$ (in any order). 
\end{enumerate}

Letting $\lambda=\omega_k=\varepsilon_1+\ldots+\varepsilon_k$, condition (R1) is easily checked; for instance, a root $(a,\overline{b})$ appears twice in $\Gamma(k)$ if $a<b\le k$, once if $a\le k<b$, and zero times otherwise. For condition (R2), we use a case by case analysis, as follows, where $a<b<c$:
\begin{enumerate}
\item $\alpha=(a,b)$, $\beta=(b,{c})$, $\gamma=(a,{c})$;
\item $\alpha=(a,b)$, $\beta=(b,\overline{c})$, $\gamma=(a,\overline{c})$;
\item $\alpha=(a,c)$, $\beta=(b,\overline{c})$, $\gamma=(a,\overline{b})$;
\item $\alpha=(b,c)$, $\beta=(a,\overline{c})$, $\gamma=(a,\overline{b})$;
\item $\alpha=(a,b)$, $\beta=(b,\overline{b})$, $\gamma=(a,\overline{a})$;
\item $\alpha=(a,\overline{a})$, $\beta=(b,\overline{b})$, $\gamma=(a,\overline{b})$.
\end{enumerate}
Case (1) is the same as in type $A$. Each of the cases (2)-(4) has the following three subcases: $k\ge c$, $b\le k<c$, and $a\le k<b$, while each of the cases (5)-(6) has the following two subcases: $k\ge b$, and $a\le k<b$. For instance, if $b\le k<c$ in Case (3), then the subsequence of $\Gamma(k)$ consisting of $\alpha,\,\beta,\,\gamma$ is $((a,\overline{b}),\,(a,c),\,(a,\overline{b}),\,(b,\overline{c}))$. 
\end{proof}

Hence, we can construct a $\lambda$-chain as a concatenation $\Gamma:=\Gamma^{\lambda_1}\ldots\Gamma^1$, where 
\begin{equation}\label{fact}\Gamma^i=\Gamma(\lambda'_i)=\Gamma_{i2}'\ldots\Gamma_{i,\lambda_i'}'\Gamma_{i1}\ldots\Gamma_{i,\lambda_i'}\,,\;\;\;\;\mbox{and }\:\Gamma_{ij}=\Gamma_j(\lambda_i')\,,\;\:\Gamma_{ij}'=\Gamma_j'\,.\end{equation}
This $\lambda$-chain is fixed for the remainder of this paper. Thus, we can replace the notation ${\mathcal A}(\Gamma)$ with ${\mathcal A}(\lambda)$.

\begin{example}\label{ex21} {\rm Consider $n=3$ and $\lambda =(3,2,1)$, for which we have the $\lambda$-chain below. The factorization of $\Gamma$ into subchains is indicated with vertical bars, while the double vertical bars separate the subchains corresponding to different columns. The underlined pairs are only relevant in Example \ref{ex21c} below.
\begin{align}\label{exlchain}\Gamma&=\Gamma_{31}\:||\:\Gamma_{22}'\Gamma_{21}\Gamma_{22}\:||\:\Gamma_{12}'\Gamma_{13}'\Gamma_{11}\Gamma_{12}\Gamma_{13}=\\&=((1,\overline{2}), \underline{(1,\overline{3})},(1,\overline{1}),(1,3),(1,2)\:||\nonumber\\ &\;\;\;\;\;\;\underline{(1,\overline{2})}\:|\:(1,\overline{3}),(1,\overline{1}),(1,3)\:|\:(1,\overline{2}),(2,\overline{3}),\underline{(2,\overline{2})},\underline{(2,3)}\:|| \nonumber\\&\;\;\;\;\;\;(1,\overline{2})\:|\:(1,\overline{3}),(2,\overline{3})\:|\:(1,\overline{1})\:|\:(1,\overline{2}), (2,\overline{2})\:|\:(1,\overline{3}),(2,\overline{3}), (3,\overline{3}))\,.\nonumber\end{align}
We represent the Young diagram of $\lambda$ inside a broken $3\times 2$ rectangle, as below. In this way, a reflection in $\Gamma$ can be viewed as swapping entries and/or  changing signs in the two parts of each column, or only the top part. 
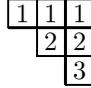
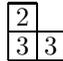
\begin{equation*}
 \begin{array}{l} \tableau{{1}&{1}&{1}\\ &{2}&{2}\\&&{3}}\\ \\
\tableau{{2}\\ {3}&{3}} \end{array}
\end{equation*}}
\end{example}

Given the $\lambda$-chain $\Gamma$ above, in Section \ref{alcovewalks} we considered subsets $J=\{ j_1<\ldots< j_s\}$ of $[m]$, where  $m$ is the length of $\Gamma$. Instead of $J$, it is now convenient to use the subsequence $T$ of the roots in $\Gamma$ whose positions are in $J$. This is viewed as a concatenation with distinguished factors $T_{ij}$ and $T_{ik}'$ induced by the factorization (\ref{fact}) of $\Gamma$. 

All the notions defined in terms of $J$ are now redefined in terms of $T$. As such, from now on we will write  $\phi(T)$, $\mu(T)$, and $|T|$, the latter being the size of $T$, cf. (\ref{defphimu}). If $J$ is a $w$-admissible subset for some $w$ in $B_n$, we will also call the corresponding $T$ a {\em $w$-admissible sequence}, and $(w,T)$ an admissible pair. We will use the notation ${\mathcal A}(\Gamma)$ and ${\mathcal A}(\lambda)$ accordingly.  We denote by $wT_{\lambda_1,1}\ldots T_{ij}$ and $wT_{\lambda_1,1}\ldots T_{ik}'$ the permutations obtained from $w$ via right multiplication by the transpositions in $T_{\lambda_1,1},\ldots, T_{ij}$ and $T_{\lambda_1,1},\ldots, T_{ik}'$, considered from left to right. This agrees with the above convention of using pairs to denote both roots and the corresponding reflections. As such, $\phi(T)$ can now be written simply $T$.

\begin{example}\label{ex21c}{\rm We continue Example \ref{ex21}, by picking the admissible pair $(w,J)$ with $w=\overline{1}\,\overline{2}\,\overline{3}\in B_3$ and $J=\{2,6,12,13\}$ (see the underlined positions in (\ref{exlchain})). Thus, we have 
\[T=T_{31}\:||\:T_{22}'T_{21}T_{22}\:||\:T_{12}'T_{13}'T_{11}T_{12}T_{13}=((1,\overline{3})\:||\:(1,\overline{2})\:|\;\;\;|\:(2, \overline{2}),(2,3)\:||\;\;\;|\;\;\;|\;\;\;|\;\;\;|\;\;\;)\,.\]
The corresponding decreasing chain in Bruhat order is the following, where the swapped entries are shown in bold (we represent permutations as broken columns starting with $w$, as discussed in Example \ref{ex21}, and we display the splitting of the chain into subchains induced by the above splitting of $T$):
\[\begin{array}{l}\tableau{{\mathbf{\overline{1}}}}\\ \\ \tableau{{\overline{2}}\\{\mathbf{\overline{3}}}} \end{array}\:>\:
\begin{array}{l} \tableau{{3}}\\ \\ \tableau{{\overline{2}}\\{1}} \end{array}\:||\:
\begin{array}{l} \tableau{{\mathbf{3}}\\{\mathbf{\overline{2}}}}\\ \\ \tableau{{1}} \end{array}\:>\:
\begin{array}{l} \tableau{{2}\\{\overline{3}}}\\ \\ \tableau{{1}} \end{array}\:|\:
\begin{array}{l} \tableau{{2}\\{\overline{3}}}\\ \\ \tableau{{1}} \end{array}\:|\:
\begin{array}{l} \tableau{{2}\\{\mathbf{\overline{3}}}}\\ \\ \tableau{{1}} \end{array}\:>\:
\begin{array}{l} \tableau{{2}\\{\mathbf{3}}}\\ \\ \tableau{{\mathbf{1}}} \end{array}\:>\:
\begin{array}{l} \tableau{{2}\\{1}}\\ \\ \tableau{{3}} \end{array}\:||\:
\begin{array}{l} \tableau{{2}\\{1}\\{3}}\\ \\ \end{array} \:|\:
\begin{array}{l} \tableau{{2}\\{1}\\{3}}\\ \\ \end{array} \:|\:
\begin{array}{l} \tableau{{2}\\{1}\\{3}}\\ \\ \end{array} \:|\:
\begin{array}{l} \tableau{{2}\\{1}\\{3}}\\ \\ \end{array} \:|\:
\begin{array}{l} \tableau{{2}\\{1}\\{3}}\\ \\ \end{array}
.
\]}
\end{example}

Given a (not necessarily admissible) pair $(w,T)$, with $T$ split into factors $T_{ij}$ and $T_{ik}'$ as above, we consider the permutations
\[\pi_{ij}=\pi_{ij}(w,T):=wT_{\lambda_1,1}\ldots T_{i,j-1}\,,\;\;\;\;\;\pi_{ik}'=\pi_{ik}'(w,T):=wT_{\lambda_1,1}\ldots T_{i,k-1}'\,;\]
when undefined, $T_{i,j-1}$ and $T_{i,k-1}'$ are given by conventions similar to (\ref{conv}), based on the corresponding factorization (\ref{fact}) of the $\lambda$-chain $\Gamma$. In particular, $\pi_{\lambda_1,1}=w$. 

\begin{definition}\label{deffill}
The {\em filling map} is the map $f$ from pairs $(w,T)$, not necessarily admissible, to fillings $\sigma=f(w,T)$ of the shape $\widehat{\lambda}$, defined (based on the notation {\em (\ref{filling})}) by
\begin{equation}\label{defswt}C_{ij}=\pi_{ij}[1,\lambda_i']\,,\;\;\;\;\;C_{ik}'=\pi_{ik}'[1,\lambda_i']\,.\end{equation}
\end{definition}

\begin{example} {\rm Given $(w,T)$ as in Example \ref{ex21c}, we have
\[f(w,T)=\tableau{{\overline{1}}&{3}&{2}&{2}&{2}&{2}&{2}\\&{\overline{2}}&{\overline{3}}&{\overline{3}}&{1}&{1}&{1}\\&&&&{3}&{3}&{3}}\,.\]}
\end{example}

The following theorem describes the way in which our tableau formula (\ref{newform}) is obtained by compressing Ram's version of Schwer's formula (\ref{hlpform}).

\begin{theorem}\label{ccompress}
{\rm (1)} We have $f(\mathcal{A}(\lambda))=\F(\lambda)$. 

{\rm (2)} Given any $\sigma \in{\mathcal F}(\lambda)$ and $(w,T)\in f^{-1}(\sigma)$, we have ${\rm ct}(f(w,T))=w(\mu(T))$. 

{\rm (3)} The following compression formula holds:
\begin{equation}\label{fcomp}\sum_{(w,T)\in f^{-1}(\sigma)}t^{\frac{1}{2}(\ell(w)+\ell(wT)-|T|)}\,(1-t)^{|T|}=t^{N(\sigma)}\,(1-t)^{\des(\sigma)}\,.\end{equation}
\end{theorem}

\begin{proof}
We start with part (1). The fact that $f(\mathcal{A}(\lambda))\subseteq\F(\lambda)$ is clear from the definition of the set of fillings ${\mathcal F}(\lambda)$ in Section \ref{sectform} and the construction (\ref{fact}) of the fixed $\lambda$-chain $\Gamma$. Viceversa, given a filling $\sigma$ in ${\mathcal F}(\lambda)$, it is not hard to construct an admissible pair $(w,T)$ in $f^{-1}(\sigma)$. We will assign to the columns $C_{ij}$ and $C_{ij}'$ signed permutations $\rho_{ij}$ and $\rho_{ij}'$ in $B_n$ recursively, starting with $\rho_{11}:=C_{11}$; in parallel, we construct the reverse ${\rm rev}(T)$ of the mentioned chain of roots $T$, and conclude by letting $w:=\rho_{\lambda_1,1}$. Each time we pass to the left neighbor $C_{ik}'$ of a column $C_{i,k+1}'=C_{ik}'(r_1,\overline{k})\ldots(r_p,\overline{k})$, we append to the part of ${\rm rev}(T)$ already constructed the roots $(r_p,\overline{k}),\ldots,(r_1,\overline{k})$ and let $\rho_{ik}':=\rho_{i,k+1}'(r_p,\overline{k})\ldots(r_1,\overline{k})$. We proceed similarly when passing to the left neighbor $C_{ij}$ of a column $C_{i,j+1}$, where $C_{i,j+1}$ differs from $D=C_{ij}(r_1,\overline{\jmath})\ldots (r_p,\overline{\jmath})$ only in position $j$; the only difference is that, in this case, we start by applying to $\rho_{i,j+1}$ and appending to $\rev(T)$ the reflection which exchanges the entry $C_{i,j+1}(j)$ with $D(j)$, and then we proceed as above.

Parts (2) and (3) of the theorem are proved in Section \ref{part2} and Section \ref{part3}, respectively.
\end{proof}

\begin{remarks}\label{specialfill} (1) The Kashiwara-Nakashima tableaux \cite{kancgr} of shape $\lambda$ index the basis elements of the irreducible representation of $\mathfrak{sp}_{2n}$ of highest weight $\lambda$. It is shown in Proposition \ref{kasnak} below that these tableaux correspond precisely to the surviving fillings in our formula (\ref{newform}) when we set $t=0$. 

(2) In (\ref{fcomp}), in general, we cannot replace the filling map $f$ with the map $\overline{f}$, sending $(w,T)$ to the compressed version $\overline{f(w,T)}$ of $f(w,T)$. Indeed, consider $n=2$, $\lambda=(3,2)$, and the following filling of $2\lambda=(6,4)$, which happens to be the ``doubled'' version of a Kashiwara-Nakashima tableau:
\[\overline{\sigma}=\tableau{{\overline{2}}&{\overline{2}}&{\overline{2}}&{\overline{2}}&{1}&{1}\\&&{\overline{1}}&{\overline{1}}&{2}&{2}}\,.\]
If $(w,T)\in\overline{f}^{-1}(\overline{\sigma})$, we must have $w=\overline{2}\overline{1}$ and 
\[T\subseteq\Gamma_{21}\Gamma_{22}=((1,\overline{1})\:|\: (1,\overline{2}), (2,\overline{2}))\,,\]
where the full $\lambda$-chain factors as follows:
\[\Gamma=\Gamma_{31}\:||\:\Gamma_{22}'\Gamma_{21}\Gamma_{22}\:||\:\Gamma_{12}'\Gamma_{11}\Gamma_{12}\,.\]
There are two elements $(w,T^1)$ and $(w,T^2)$ in $\overline{f}^{-1}(\overline{\sigma})$, namely
\[T^1=((1,\overline{2}))\,,\;\;\;\; \mbox{and}\;\;\;\; T^2=((1,\overline{1}), (1,\overline{2}), (2,\overline{2}))\,.\]
But we have
\[\sum_{(w,T)\in \overline{f}^{-1}(\overline{\sigma})}t^{\frac{1}{2}(\ell(w)+\ell(wT)-|T|)}\,(1-t)^{|T|}=t(1-t)+(1-t)^3=(1-t)(1-t+t^2)\,.\]
In general, the above sum has several factors not of the form $t$ or $1-t$, which are hard to control.

(3) In order to measure the compression phenomenon,  we define the {\em compression factor} $c(\lambda)$ like in \cite{lenhlp}, as the ratio of the number of terms in Ram's version of Schwer's formula for $\lambda$ and the number of terms in the tableau formula. The compression factor is considerably larger in type $C$. For instance, for $\lambda=(3, 2, 1, 0)$ in $C_4$ we have 23,495 terms in the compressed formula, while $c(\lambda)=44.9$. 
\end{remarks}

\begin{proposition}\label{kasnak} The map $\sigma\mapsto \overline{\sigma}$ (see {\rm (\ref{redfilling})}) is a bijection between the fillings $\sigma$ in ${\mathcal F}(\lambda)$ with $N(\sigma)=0$ and the ``doubled'' versions of the type $C$ Kashiwara-Nakashima tableaux of shape $\lambda$.
\end{proposition}

\begin{proof}
It was proved in \cite{laaapm} that for each type $C$ Kashiwara-Nakashima tableau of shape $\lambda$ there is a unique admissible pair $(w,T)$ whose associated chain in Bruhat order is saturated and ends at the identity, such that the compressed version $\overline{\sigma}$ of $\sigma=f(w,T)$ is the ``doubled'' version of the given tableau. It follows that the term associated to $(w,T)$ in (\ref{hlpform}) is $t^0(1-t)^{|T|}x^{w(\mu(T))}$, and therefore $N(\sigma)=0$, by (\ref{fcomp}). On the other hand, since $P_\lambda(x;0)$ is the irreducible character indexed by $\lambda$, which is expressed in terms of Kashiwara-Nakashima tableaux, we conclude that all $\sigma$ in ${\mathcal F}(\lambda)$ with $N(\sigma)=0$ arise in this way.
\end{proof}

\section{The tableau formula in type $B$}\label{typeb} We now restrict ourselves to the root system of type $B_n$. 
This can be represented as 
$\Phi=\{\pm\varepsilon_i\pm\varepsilon_j \::\:  1\leq i<j\leq n\}\cup\{\pm \varepsilon_i\::\: 1\leq i\leq n\}$. 
The simple roots are $\alpha_i=\varepsilon_i-\varepsilon_{i+1}$, 
for $i=1,\ldots,n-1$ and $\alpha_n=\varepsilon_n$. 
The fundamental weights are $\omega_i = \varepsilon_1+\ldots +\varepsilon_i$, 
for $i=1,\ldots,n-1$ and $\omega_n=\frac12 (\varepsilon_1+\cdots+\varepsilon_n)$. A dominant weight $\lambda=\alpha_1\omega_1+\ldots+\alpha_n\omega_n$, where $\alpha_i\in\Z_{\ge 0}$, is identified with the partition $\mu=(n^{\alpha_n},\ldots,1^{\alpha_1})$; we let $\ell(\mu):=\alpha_1+\ldots+\alpha_n$, and write $\mu=(\mu_1,\ldots,\mu_{\ell(\mu)})$. A dominant weight is regular if $\alpha_i>0$ for all $i$. Let us now fix such a weight $\lambda$.

The corresponding Weyl group $W$ is the same group of signed permutations $B_n$ considered above. For simplicity, we again use the same notation for roots and the corresponding reflections, cf. Section \ref{sectform}. The pairs $(i,j)$ and $(i,\overline{\jmath})$ have the same meaning as in type $C$, whereas $(i)$ denotes the positive root $\varepsilon_i$. Note that, as a reflection in $B_n$, $(i)$ is the same as $(i,\overline{\imath})$ in type $C$.

The canonical $\omega_k$-chains and $\lambda$-chains are very similar to those in type $C$. If $k<n$, let 
\[\Gamma(k):=\Gamma_1'\ldots\Gamma_k'\Gamma_1(k)\ldots\Gamma_k(k)\,,\]
where
\begin{align*}
&\;\:\Gamma_j':=((1,\overline{\jmath}),(2,\overline{\jmath}),\ldots,(j-1,\overline{\jmath}),(j))\,,\\
&\begin{array}{llllll}\Gamma_j(k):=(\!\!\!\!\!&(1,\overline{\jmath}),&(2,\overline{\jmath}),&\ldots,&(j-1,\overline{\jmath}),\\
&(j,\overline{k+1}),&(j,\overline{k+2}),&\ldots,&(j,\overline{n}),\\&(j),\\ 
&(j,n),&(j,n-1),&\ldots,&(j,k+1)\,)\,.\end{array}\end{align*}
On the other hand, we let 
\[\Gamma(n):=\Gamma_1'\ldots\Gamma_n'=\Gamma_1(n)\ldots\Gamma_n(n)\,.\]
Like in the type $C$ case, we can prove that $\Gamma(k)$ is an $\omega_k$-chain for any $k$. Hence, we can construct a $\lambda$-chain as a concatenation $\Gamma:=\Gamma^{\ell(\mu)}\ldots\Gamma^1$, where $\Gamma^i=\Gamma(\mu_i)$.

The filling map is defined like in Definition \ref{deffill}. This gives rise to fillings 
\begin{equation*}\sigma={\mathcal C}^{\ell(\mu)}\ldots {\mathcal C}^{1}\,,\end{equation*}
where each ${\mathcal C}^{i}$ is a concatenation of columns of height $\mu_i$, as follows:
\[{\mathcal C}^{i}:=\casethree{C_{i1}'\ldots C_{i,\mu_i}'C_{i1}\ldots C_{i,\mu_i}}{\mu_i<n}{C_{i1}\ldots C_{i,\mu_i}}{i\ne 1\mbox{ and }\mu_i=n}{C_{11}}{i=1}\]
The fillings are subject to the same conditions (1)-(3) as in type $C$, where condition (3) is made more precise below. In fact, the above $\lambda$-chain $\Gamma$ specifies the way in which each column is related to its left neighbor. Essentially, everything is similar to type $C$, except for a small difference in the passage from some column $C_{ik}'$ to $C_{i,k+1}'$. Namely, there exist positions $1\le r_1<\ldots <r_p<k$ (possibly $p=0$) such that $C_{i,k+1}'=C_{ik}'(r_1,\overline{k})\ldots (r_p,\overline{k})$, like in type $C$, or $C_{i,k+1}'=C_{ik}'(r_1,\overline{k})\ldots (r_p,\overline{k})(k)$, in which case we also require $C_{i,k+1}'(k)\le n$. 

The weight of a filling, and the statistics $N(\sigma)$ and $\des(\sigma)$ are defined completely similarly to type $C$. The only minor addition is the definition of $N(D,C)$ and $\des(D,C)$ when $C=D(r_1,\overline{k})\ldots (r_p,\overline{k})(k)$. With the notation in Case 2 of the definition of $N(D,C)$ in Section \ref{sectform}, we set
\[N(D,C):=N(D,D_p)+N(D_p,C)\,,\;\;\;\des(D,C):=p+1\,.\]
Here $N(D,D_p)$ is defined in Case 2 above, whereas $N(D_p,C)$ is given by the second formula in (\ref{diflengths}); more precisely,
\[N(D_p,C):=\frac{1}{2}|(\overline{a},a)\setminus\{\pm D_p(i)\::\:i=1,\ldots,k\}|\,,\]
where $a:=D_p(k)$.

Given the above constructions, the proof of the following theorem is completely similar to its counterparts in type $C$, since no new situations arise.

\begin{theorem} Theorems {\rm \ref{newformc}} and {\rm \ref{ccompress}} hold in type $B$, with the appropriate constructions explained above.
\end{theorem}

\begin{remark}
The situation in type $D$ is slightly more complex. In this case, applying the above ideas leads to an analog of the compression formula (\ref{fcomp}) which contains factors of the form $1-t^k$ with $k>1$ in the right-hand side. However, these factors are not hard to control, while no extra factors appear. 
\end{remark}

\section{Proof of Theorem {\ref{ccompress} (2)}}\label{part2}

Recall the $\lambda$-chain $\Gamma$ in Section \ref{specschwer}. Let us write $\Gamma=(\beta_1,\ldots,\beta_m)$, as in Section \ref{alcovewalks}. As such, we recall the hyperplanes $H_{\beta_k,l_k}$ and the corresponding affine reflections $\widehat{r}_k=s_{\beta_k,l_k}=s_{\beta_k}+l_k\beta_k$. 

Now fix a signed permutation $w$ in $B_n$ and a subset $J=\{j_1<\ldots<j_s\}$ of $[m]$ (not necessarily $w$-admissible).  Let $\Pi$ be the alcove path corresponding to $\Gamma$, and define the alcove walk $\Omega$ as in Section \ref{alcovewalks}, by
\[\Omega:=\phi_{j_1}\ldots \phi_{j_s} (w(\Pi))\,.\]
Given $k$ in $[m]$, let $i=i(k)$ be the largest index in $[s]$ for which $j_i<k$, and let $\gamma_k:=wr_{j_1}\ldots r_{j_i}(\beta_k)$. Then the hyperplane containing the face $F_k$ of $\Omega$ (cf. Definition \ref{defalcwalk}) is of the form $H_{\gamma_k,m_k}$; in other words
\[H_{\gamma_k,m_k}=w\widehat{r}_{j_1}\ldots \widehat{r}_{j_i}(H_{\beta_k,l_k})\,.\] 
Our first goal is to describe $m_k$ purely in terms of the filling associated to $(w,J)$.

Let $\widehat{t}_k$ be the affine reflection in the hyperplane $H_{\gamma_k,m_k}$. Note that
\[\widehat{t}_k=w\widehat{r}_{j_1}\ldots \widehat{r}_{j_i}\widehat{r}_k\widehat{r}_{j_i}\ldots \widehat{r}_{j_1}w^{-1}\,.\]
Thus, we can see that
\[w\widehat{r}_{j_1}\ldots \widehat{r}_{j_i}=\widehat{t}_{j_i}\ldots \widehat{t}_{j_1}w\,.\]

Let $T=((a_1,b_1),\ldots,(a_s,b_s))$ be the subsequence of $\Gamma$ indexed by the positions in $J$, cf. Section \ref{specschwer}. Let $T^i$ be the initial segment of $T$ with length $i$, let $w_i:=wT^i$, and $\sigma_i:=\overline{f(w,T^i)}$, see (\ref{redfilling}). In particular, $\sigma_0$ is the filling with all entries in row $i$ equal to $w(i)$, and  $\sigma:=\sigma_s=\overline{f(w,T)}$. The columns of a filling of $2\lambda$ are numbered left to right by $2\lambda_1$ to 1. We split each segment $\Gamma^k$ of $\Gamma$ into two parts: the head which is a concatenation of $\Gamma'_{k\cdot}$, and the tail which is a concatenation of $\Gamma_{k\cdot}$, see (\ref{fact}). We say that the head corresponds to column $2k$ of the Young diagram $2\lambda$, whereas the tail corresponds to column $2k-1$ (cf. the construction of ${f(w,T)}$ in Section \ref{specschwer} and (\ref{redfilling})). If $\beta_{j_{i+1}}=(a_{i+1},b_{i+1})=(a,b)$ falls in the segment of $\Gamma$ corresponding to column $p$ of $2\lambda$, then $\sigma_{i+1}$ is obtained from $\sigma_i$ by replacing the entry $w_i(a)$ with $w_i(b)$ in the columns $p-1,\ldots,1$ of $\sigma_i$, as well as, possibly, the entry $w_i(\overline{b})$ with ${w_i(\overline{a})}$ in the same columns. 

Now fix a position $k$, and consider $i=i(k)$ and the roots $\beta_k$, $\gamma:=\gamma_k$, as above, where $\gamma_k$ might be negative. Assume that $\beta_k$ falls in the segment of $\Gamma$ corresponding to column $q$ of $2\lambda$. Given a filling $\phi$, we denote by $\phi[p]$ and $\phi(p,q]$ the parts of $\phi$ consisting of the columns $2\lambda_1,2\lambda_1-1,\ldots,p$ and $p-1,p-2,\ldots,q$, respectively. We also recall the definition (\ref{defct}) and conventions related to the content of a filling; this definition now applies to any filling of some Young diagram. 

\begin{proposition}\label{level} With the above notation, we have
\[m_k=\langle{\rm ct}(\sigma[q]),\gamma^\vee\rangle\,.\]
\end{proposition}

\begin{proof}  We apply induction on $i$, which starts at $i=0$, when the verification is straightforward. We will now proceed from $j_1<\ldots<j_{i}<k$, where $i=s$ or $k\le j_{i+1}$, to $j_1<\ldots<j_{i+1}<k$, and we will freely use the notation above. 
Assume that $\beta_{j_{i+1}}$ falls in the segment of $\Gamma$ corresponding to column $p$ of $2\lambda$, where $p\ge q$. 

We need to compute
\[w\widehat{r}_{j_1}\ldots \widehat{r}_{j_{i+1}}(H_{\beta_k,l_k})=\widehat{t}_{j_{i+1}}\ldots \widehat{t}_{j_1}w(H_{\beta_k,l_k})=\widehat{t}_{j_{i+1}}(H_{\gamma,m})\,,\]
where $m=\langle{\rm ct}(\sigma_i[q]),\gamma^\vee\rangle$, by induction. Let $\gamma':=\gamma_{j_{i+1}}$, and $\widehat{t}_{j_{i+1}}=s_{\gamma',m'}$, where $m'=\langle{\rm ct}(\sigma_i[p]),(\gamma')^\vee\rangle$, by induction. We will use the following formula:
\[s_{\gamma',m'}(H_{\gamma,m})=H_{s_{\gamma'}(\gamma),m-m'\langle\gamma',\gamma^\vee\rangle}\,.\]
Thus, the proof is reduced to showing that
\[m-m'\langle\gamma',\gamma^\vee\rangle=\langle{\rm ct}(\sigma_{i+1}[q]),s_{\gamma'}(\gamma^\vee)\rangle\,.\] 
An easy calculation, based on the above information, shows that the latter equality is non-trivial only if $p>q$, in which case it is equivalent to
\begin{equation*}\langle{\rm ct}(\sigma_{i+1}(p,q])-{\rm ct}(\sigma_{i}(p,q]),\gamma^\vee\rangle=\langle\gamma',\gamma^\vee\rangle\,\langle{\rm ct}(\sigma_{i+1}(p,q]),(\gamma')^\vee\rangle\,.\end{equation*}
This equality is a consequence of the fact that
\[{\rm ct}(\sigma_{i+1}(p,q])=s_{\gamma'}({\rm ct}(\sigma_{i}(p,q]))\,,\]
which follows from the construction of $\sigma_{i+1}$ from $\sigma_i$ explained above.
\end{proof}

\begin{proof}[Proof of Theorem {\rm {\ref{ccompress} (2)}}] We apply induction on the size of $T$, using freely the notation above. We prove the statement for $T=(\beta_{j_1},\ldots,\beta_{j_{s+1}})$, assuming it holds for $T^s=(\beta_{j_1},\ldots,\beta_{j_{s}})$. We have
\[w(\mu(T)=w\widehat{r}_{j_1}\ldots \widehat{r}_{j_{s+1}}(\lambda)=\widehat{t}_{j_{s+1}}\ldots \widehat{t}_{j_1}w(\lambda)=\widehat{t}_{j_{s+1}}({\rm ct}(\sigma_s))\,,\]
by induction. We need to show that
\begin{equation}\label{toprove}\widehat{t}_{j_{s+1}}({\rm ct}(\sigma_s))={\rm ct}(\sigma_{s+1})\,.\end{equation}
Let $\gamma:=\gamma_{j_{s+1}}$ and assume that $\beta_{j_{s+1}}$ falls in the segment of $\Gamma$ corresponding to column $p$ of $2\lambda$. Based on Proposition \ref{level}, (\ref{toprove}) is rewritten as
\begin{equation}\label{toprove1}s_\gamma({\rm ct}(\sigma_s))+\langle{\rm ct}(\sigma_s[p]),\gamma^\vee\rangle\gamma={\rm ct}(\sigma_{s+1})\,.\end{equation}
Decomposing ${\rm ct}(\sigma_s)$ as ${\rm ct}(\sigma_s[p])+{\rm ct}(\sigma_s(p,1])$ (cf. the notation above), and ${\rm ct}(\sigma_{s+1})$ in a similar way, (\ref{toprove1}) would follow from
\begin{align*}&s_\gamma({\rm ct}(\sigma_s[p]))+\langle{\rm ct}(\sigma_s[p]),\gamma^\vee\rangle\gamma={\rm ct}(\sigma_{s+1}[p])\,,\\&s_\gamma({\rm ct}(\sigma_s(p,1]))={\rm ct}(\sigma_{s+1}(p,1])\,.\end{align*}
The first equality is clear since $\sigma_s[p]=\sigma_{s+1}[p]$, while the second one follows from the construction of $\sigma_{s+1}$ from $\sigma_s$ explained above.
\end{proof}

\section{Proof of Theorem {\ref{ccompress} (3)}}\label{part3}

We start by recalling some basic facts about the group $B_n$ and some notation from Section \ref{sectform}. We will use the following notation related to a word $w=w_1\ldots w_l$ with integer letters, which is sometimes identified with its set of letters:
\[w[i,j]:=w_i\ldots w_j\,,\;\;\;N_{ab}(w):=|(a,b)\cap w|\,,\;\;\;N_{ab}(\pm w):=N_{ab}(w)+N_{ab}(-w)\,,\]
where $-w:=\overline{w_1}\ldots \overline{w_l}$. Given words $w^1,\ldots,w^p$, we let
\[N_{ab}(w^1,\ldots,w^p):=N_{ab}(w^1)+\ldots+N_{ab}(w^p)\,.\]
We also let
\[\tau_{ab}:=\casetwo{1}{a,b\le n}{0}\]
With this notation, given a signed permutation $w$ in $B_n$ and $1\le i<j\le n$, $a:=w(i)$, $b:=w(j)$, we have the following facts:
\begin{align}\label{diflengths}
&\frac{\ell(w(i,j))-\ell(w)-1}{2}=N_{ab}(w[i,j])\,,\nonumber\\
&\frac{\ell(w(i,\overline{\imath}))-\ell(w)-1}{2}=N_{a\overline{a}}(w[i,n])\,,\\
&\frac{\ell(w(i,\overline{\jmath}))-\ell(w)-1}{2}=N_{a\overline{b}}(w[i,j-1],\pm w[j+1,n])+\tau_{ab}\,,\nonumber
\end{align}
assuming that the left-hand side is nonnegative (i.e., that we go up in Bruhat order upon applying the indicated reflection); these facts are used implicitly throughout.

Given a chain of roots $\Delta$, we define ${\mathcal A}^r(\Delta)$ like in (\ref{decch}) except that we impose an increasing chain condition and $w\in W$. The following simple lemma will be useful throughout, for splitting chains into subchains.

\begin{lemma}\label{split}
Consider $(w,T)$ with $T$ written as a concatenation $S_1\ldots S_p$. Let $w_i:=wS_1\ldots S_i$, so $w_0=w$. Then
\[\frac{1}{2}(\ell(w)+\ell(wT)-|T|)=\frac{1}{2}(\ell(w_{p-1})+\ell(w_p)-|S_p|)+\sum_{i=1}^{p-1}\frac{1}{2}(\ell(w_{i-1})-\ell(w_i)-|S_i|)\,.\]
\end{lemma}

Let $\Delta$ be the chain
\[
\begin{array}{llllll}\Delta:=(\!\!\!\!\!&(1,p+1),&(1,p+2),&\ldots,&(1,n),\\&(1,\overline{1}),\\
&(1,\overline{n}),&(1,\overline{n-1}),&\ldots,&(1,\overline{p+1})\,)\,.\end{array}
\]

\begin{proposition}\label{p2cols0}
Consider a signed permutation $w$ in $B_n$ with $a:=w(1)$, a position $1\le p\le n$, and a value $b\in\{\pm a\}\cup (\pm w[p+1,n])$ such that $b\ge a$. Then we have
\begin{equation}\label{sum2cols0}
\stacksum{T\::\:(w,T)\in{\mathcal A}^r(\Delta)}{wT(1)=b} t^{\frac{1}{2}(\ell(wT)-\ell(w)-|T|)} (1-t)^{|T|}=t^{N_{ab}(w[2,p])} (1-t)^{1-\delta_{ab}}\,;
\end{equation}
here $\delta_{ab}$ is the Kronecker delta.
\end{proposition}

\begin{proof}
 Given $s\in\{\overline{1},\pm(p+1),\ldots,\pm n\}$, we let $\Delta_s$ be the subchain of $\Delta$ starting with $(1,s)$. We also let
\[S(w,s):=\stacksum{T\::\:(w,T)\in{\mathcal A}^r(\Delta_{s})}{wT(1)=b} t^{\frac{1}{2}(\ell(wT)-\ell(w)-|T|)} (1-t)^{|T|}\,.\]
We consider three cases, depending on  $b=w(q)$, $b=\overline{w(q)}$, and $b=\overline{a}$.
The proof in the first case is identical to that of the analogous result for type $A$, namely \cite{lenhlp}[Proposition 5.3]. 

{\em Case {\rm 2}.} Let $c:=w(q)=\overline{b}$, and $p<q\le s$. We start by showing
\begin{equation}\label{ind1}S(w,\overline{s})=t^{N_{a\overline{c}}(w[2,q-1],w[q+1,s],\pm w[s+1,n])+\tau_{ac}} (1-t)\,.\end{equation}
We use induction on $s$, which starts at $s=q$. For $s>q$, let $w^1:=w[2,q-1]$, $w^2:=w[q+1,s-1]$, $w^3:=w[s+1,n]$, and $d:=w(s)$. The sum ${S}(w,\overline{s})$ splits into two sums, depending on $(1,\overline{s})\not\in T$, and $(1,\overline{s})\in T$. By induction, the first sum is
\[{S}(w,\overline{s-1})=t^{N_{a\overline{c}}(w^1,w^2,\pm dw^3)+\tau_{ac}}(1-t)\,.\]
Again by induction, if $a<\overline{d}<\overline{c}$, then the second sum is
\begin{align*}&t^{N_{a\overline{d}}(w^1cw^2,\pm w^3)+\tau_{ad}}(1-t){S}(w(1,\overline{s}),\overline{s-1})\\=&t^{N_{a\overline{d}}(w^1cw^2,\pm w^3)+N_{\overline{d}\overline{c}}(w^1,w^2,\pm \overline{a}w^3)+\tau_{ad}+\tau_{\overline{d}c}}(1-t)^2\,;\end{align*}
otherwise, it is empty. Adding up the two sums into which ${S}(w,\overline{s})$ splits, we obtain
\[t^{N_{a\overline{c}}(w^1,w^2d,\pm w^3)+\tau_{ac}}(1-t)\,.\]
The last claim rests on the easily verified facts that, if $a<\overline{d}<\overline{c}$, then
\[\tau_{ad}+\tau_{\overline{d}c}=\tau_{ac}\,,\;\;\;\;\; N_{a\overline{d}}(c)+N_{\overline{d}\overline{c}}(\overline{a})=N_{a\overline{c}}(d)\,.\]

Still assuming that $c=w(q)=\overline{b}$ and $p<q$, we now show that 
\begin{equation}\label{ind2}S(w,\overline{1})=t^{N_{a\overline{c}}(w[2,q-1],w[q+1,n])+\tau'_{ac}} (1-t)\,,\end{equation}
where 
\[\tau'_{ac}:=\casetwo{1}{a<c\le n}{0}\]
Let $w^1:=w[2,q-1]$, like before, and $w^2:=w[q+1,n]$. The sum ${S}(w,\overline{1})$ splits into two sums, depending on $(1,\overline{1})\not\in T$, and $(1,\overline{1})\in T$. By (\ref{ind1}), the first sum is
\[{S}(w,\overline{n})=t^{N_{a\overline{c}}(w^1,w^2)+\tau_{ac}}(1-t)\,.\]
Again by (\ref{ind1}), if $c<a\le n$, then the second sum is
\[t^{N_{a\overline{a}}(w^1cw^2)}(1-t){S}(w(1,\overline{1}),\overline{n})=t^{N_{a\overline{a}}(w^1cw^2)+N_{\overline{a}\overline{c}}(w^1,w^2)+\tau_{\overline{a}c}}(1-t)^2\,;\]
otherwise, it is empty. Adding up the two sums into which ${S}(w,\overline{s})$ splits, we obtain
\[t^{N_{a\overline{c}}(w^1,w^2)+\tau'_{ac}}(1-t)\,.\]

Assuming that $c=w(q)=\overline{b}$ and $p<q<s$, we now show that 
\begin{equation}\label{ind3}S(w,s)=t^{N_{a\overline{c}}(w[2,q-1],w[q+1,s-1])+\tau'_{ac}} (1-t)\,.\end{equation}
We use decreasing induction on $s$. Like before, we let $w^1:=w[2,q-1]$, $w^2:=w[q+1,s-1]$, and $d:=w(s)$. The sum ${S}(w,{s})$ splits into two sums, depending on $(1,{s})\not\in T$, and $(1,{s})\in T$. By induction, the first sum is
\[{S}(w,{s+1})=t^{N_{a\overline{c}}(w^1,w^2d)+\tau'_{ac}}(1-t)\,.\]
Again by induction, if $a<{d}<\overline{c}$, then the second sum is
\[t^{N_{a{d}}(w^1cw^2)}(1-t){S}(w(1,{s}),{s+1})=t^{N_{a{d}}(w^1cw^2)+N_{{d}\overline{c}}(w^1,w^2a)+\tau'_{{d}c}}(1-t)^2\,;\]
otherwise, it is empty. (In both calculations, induction works by substituting $\overline{1}$ for $n+1$ when $s=n$, and by using (\ref{ind2}) in this case.) Adding up the two sums into which ${S}(w,{s})$ splits, we obtain
\[t^{N_{a\overline{c}}(w^1,w^2)+\tau'_{ac}}(1-t)\,.\]
The last claim rests on the easily verified fact that, if $a<{d}<\overline{c}$, then
\[N_{a{d}}(c)+\tau'_{dc}=\tau'_{ac}\,.\]

{\em Case {\rm 3}.} Let us now assume that $b=\overline{a}$. We need to show that
\begin{equation}\label{ind41}S(w,p+1)=t^{N_{a\overline{a}}(w[2,p])}(1-t)\,.\end{equation}
We use decreasing induction on $p$, which starts at $p=n$; in this case $\Delta$ only contains the pair $(1,\overline{1})$, so the above convention of substituting $\overline{1}$ for $n+1$ works well here too. For $p<n$, we let $d:=w(p+1)$. The sum ${S}(w,p+1)$ splits into two sums, depending on $(1,p+1)\not\in T$, and $(1,p+1)\in T$. By induction, the first sum is
\[{S}(w,p+2)=t^{N_{a\overline{a}}(w[2,p]d)}(1-t)\,.\]
If $a<d<\overline{a}$, then by (\ref{ind3}) of Case 2, the second sum is
\[t^{N_{ad}(w[2,p])}(1-t){S}(w(1,p+1)),p+2)=t^{N_{ad}(w[2,p])+N_{d\overline{a}}(w[2,p])+\tau'_{da}}(1-t)^2\,;\]
otherwise, it is empty. Adding up the two sums into which ${S}(w,p+1)$ splits, we obtain the desired result.

{\em Case {\rm 2} (continued).} Assuming that $c=w(q)=\overline{b}$ and $p<q$, we now show that 
\begin{equation}\label{ind4}S(w,q)=t^{N_{a\overline{c}}(w[2,q-1])} (1-t)\,.\end{equation}
The sum ${S}(w,q)$ splits into two sums, depending on $(1,q)\not\in T$, and $(1,q)\in T$. By (\ref{ind3}) of Case 2, the first sum is
\[{S}(w,q+1)=t^{N_{a\overline{c}}(w[2,q-1])+\tau'_{ac}}(1-t)\,.\]
If $a<c\le n$, then by (\ref{ind41}) of Case 3, the second sum is
\[t^{N_{ac}(w[2,q-1])}(1-t){S}(w(1,q)),q+1)=t^{N_{ac}(w[2,q-1])+N_{c\overline{c}}(w[2,q])}(1-t)^2\,;\]
otherwise, it is empty. Adding up the two sums into which ${S}(w,q)$ splits, we obtain the desired result.

The final step in Case 2 is to prove that 
\begin{equation}\label{ind5}S(w,p+1)=t^{N_{a\overline{c}}(w[2,p])} (1-t)\,.\end{equation}
This can be done by decreasing induction on $p$, starting with $p=q-1$, which is the case proved in (\ref{ind4}). The procedure is completely similar to the ones above, and, in fact, to the one for type $A$ in \cite{lenhlp}[Proposition 5.3]. 
\end{proof}

Let us consider the chain
\begin{equation}\label{chphi}
\begin{array}{llllll}\Phi:=\Gamma_1(n)\ldots\Gamma_n(n)=(\!\!\!\!\!&(1,\overline{1}),\\&(1,\overline{2}),&(2,\overline{2}),\\
&&\ldots\\
&(1,\overline{n}),&(2,\overline{n}),&\ldots,&(n-1,\overline{n})\,)\,.\end{array}
\end{equation}
We denote by $\Phi_{ij}$ the subchain of $\Phi$ starting with $(i,\overline{\jmath})$.  
Given a signed permutation $w$, recall the definition (\ref{defl}) of $\ell_+(w)$ and $\ell_-(w)$. 
Given $(i,j)$ with $1\le i\le j\le n$, we also define
\begin{align}\label{defl3}&\ell^{ij}_-(w):=|\{(k,l)\::\:(k,\overline{l})\in\Phi\setminus\Phi_{ij}\,,\:w(k)>\overline{w(l)}\}|\,,\\&\overline{\ell}^{ij}_-:=\ell_-(w)-\ell^{ij}_-(w)\,.\nonumber\end{align}

\begin{proposition}\label{descpart} Fix $(i,j)$ with $1\le i\le j\le n$ and a signed permutation $w$ in $B_n$. We have
\begin{equation}\label{comp1}\sum_{T\::\:(w,T)\in{\mathcal A}(\Phi_{ij})}t^{\frac{1}{2}(\ell(w)+\ell(wT)-|T|)}\,(1-t)^{|T|}=t^{\ell_+(w)+\ell^{ij}_-(w)}\,.\end{equation}
In particular, if the above sum is over $(w,T)\in{\mathcal A}(\Phi)$, then the right-hand side is $t^{\ell_+(w)}$.
\end{proposition}

\begin{proof}
Let us denote the sum in the left-hand side of (\ref{comp1}) by $S(w,i,j)$, and the corresponding sum over $(w,T)\in{\mathcal A}(\Phi_{ij}\setminus\{(i,\overline{\jmath})\})$ by $S'(w,i,j)$. We view the chain $\Phi$ as a total order on the pairs $(i,\overline{\jmath})$, with $(1,\overline{1})$ being the smallest pair. With this in mind, we use decreasing induction on pairs $(i,\overline{\jmath})$. Given such a pair, if $w(i)<\overline{w(j)}$ then the induction step is clear, so assume the contrary. We can now split $S(w,i,j)$ into two sums, depending on $(i,\overline{\jmath})\not\in T$ and $(i,\overline{\jmath})\in T$. By induction, the first sum is
\[S'(w,i,j)=t^{1+\ell_+(w)+\ell^{ij}_-(w)}\,.\]
By induction and Lemma \ref{split}, the second sum is
\begin{align*}&(1-t)t^{\frac{1}{2}(\ell(w)-\ell(w(i,\overline{\jmath}))-1)}S'(w(i,\overline{\jmath}),i,j)\\=&(1-t)t^{\frac{1}{2}(\ell(w)-\ell(w(i,\overline{\jmath}))-1)+\ell_+(w(i,\overline{\jmath}))+\ell^{ij}_-(w(i,\overline{\jmath}))}\,.\end{align*}
The induction step is completed once we show that
\[\ell_+(w)+\ell^{ij}_-(w)=\frac{1}{2}\left(\ell(w)-\ell(w(i,\overline{\jmath}))-1\right)+\ell_+(w(i,\overline{\jmath}))+\ell^{ij}_-(w(i,\overline{\jmath}))\,.\]
The latter equality can be rewritten as
\[\Delta\ell_+(w)+\Delta\ell^{ij}_-(w)-1=\Delta\overline{\ell}^{ij}_-(w)\,,\]
where $\Delta\ell_+(w):=\ell_+(w)-\ell_+(w(i,\overline{\jmath}))$, and similarly for the other two variations. In order to prove this, let us first fix $k$ between $i$ and $j$, and analyze the contribution to the three variations of the pairs $(i,k)$ and $(k,j)$, cf. (\ref{defl}) and (\ref{defl3}). For simplicity, let $a:=w(i)$, $b:=w(k)$, and $c:=w(j)$, where $a>\overline{c}$. The mentioned nonzero contributions are as follows:
\begin{enumerate}
\item the pair $(i,k)$ contributes $1$ to $\Delta\ell_+(w)$ if $a>b>\overline{c}$;
\item the pair $(k,j)$ contributes $-1$ to $\Delta\ell_+(w)$ if $\overline{a}<b<{c}$, which is equivalent to $a>\overline{b}>\overline{c}$;
\item the pair $(i,k)$ contributes $1$ to $\Delta\ell^{ij}_-(w)$ if $a>\overline{b}>\overline{c}$;
\item the pair $(k,j)$ contributes $1$ to $\Delta\overline{\ell}^{ij}_-(w)$ if $a>b>\overline{c}$.
\end{enumerate}
Note that the second and third contributions cancel out, whereas the first one is equal to the fourth one. The analysis is completely similar if $k<i$ or $k>j$. The pair $(i,j)$ only contributes $1$ to $\Delta\ell^{ij}_-(w)$. As far as the pairs $(i,i)$ and $(j,j)$ are concerned, the contribution of the first one to $\Delta\ell^{ij}_-(w)$ and of the second one to $\Delta\overline{\ell}^{ij}_-(w)$ are both equal to $\sigma_{ac}$, see (\ref{defsigma}).
\end{proof}

\begin{proof}[Proof of Theorem {\rm \ref{ccompress} (3)}]
Fix a filling $\sigma$ in ${\mathcal F}(\lambda)$ with columns $C_{ij}$ and $C_{ij}'$, as explained in Section \ref{sectform}. Recall the chain $\Phi:=\Gamma_1(n)\ldots\Gamma_n(n)=\Gamma_{11}\ldots\Gamma_{1n}$ in (\ref{chphi}). By splitting the $\lambda$-chain $\Gamma$ into the tail $\Phi$ and its complement, and by using Lemma \ref{split}, the sum in the left-hand side of (\ref{fcomp}) can be rewritten as
\begin{align}\label{split2}&\sum_{(w,T)\in f^{-1}(\sigma)}t^{\frac{1}{2}(\ell(w)+\ell(wT)-|T|)}\,(1-t)^{|T|}=\\=&\left( \stacksum{(w,T)\in f^{-1}(\sigma)}{T_{11}=\ldots=T_{1n}=\emptyset}t^{\frac{1}{2}(\ell(w)-\ell(wT)-|T|)}\,(1-t)^{|T|}\right)\times\nonumber\\\times&\left( \sum_{T\::\:(C_{11},T)\in{\mathcal A}(\Phi)}t^{\frac{1}{2}(\ell(C_{11})+\ell(C_{11}T)-|T|)}\,(1-t)^{|T|}\right)\,;\nonumber\end{align}
here the column $C_{11}$, which has height $n$, is viewed as a signed permutation in $B_n$. By Proposition \ref{descpart}, the second bracket is $t^{\ell_+(C_{11})}$. 

In order to evaluate the first bracket, we will reverse all chains. Let us start by recalling the construction (\ref{fact}) of the $\lambda$-chain $\Gamma$, and in particular the order in which the subchains $\Gamma_{ij}$ and $\Gamma_{ij}'$ are concatenated (including the conventions in Section \ref{sectform} related to $\Gamma_{i,j+1}$ and $\Gamma_{i,j+1}'$). We denote by $\Gamma_{ij}^r$ and $(\Gamma_{ij}')^r$ the corresponding reverse chains. Also recall that we defined ${\mathcal A}^r(\,\cdot\,)$ like in (\ref{decch}) except that we imposed an increasing chain condition and $w\in W$. We consider pairs $(w_{ij},S_{ij})$ in ${\mathcal A}^r(\Gamma_{ij}^r)$ and $(w_{ij}',S_{ij}')$ in ${\mathcal A}^r((\Gamma_{ij}')^r)$, where $w_{ij}$ and $w_{ij}'$ are defined by
\[w_{ij}:=C_{11}S_{1,\lambda_1'}'\ldots S_{i,j+1}\,,\;\;\;\;\;w_{ij}':=C_{11}S_{1,\lambda_1'}'\ldots S_{i,j+1}\,,\]
where the concatenation order for $S_{ij}$ and $S_{ij}'$ comes from that for $\Gamma_{ij}$ and $\Gamma_{ij}'$; in particular, $w_{1,\lambda_1'}'=C_{11}$. Given this notation, we define the sum
\[\Sigma_{ij}:=\stacksum{S_{ij}\::\:(w_{ij},S_{ij})\in{\mathcal A}^r(\Gamma_{ij}^r)}{w_{ij}S_{ij}[1,\lambda_i']=C_{ij}} t^{\frac{1}{2}(\ell(w_{ij}S_{ij})-\ell(w_{ij})-|S_{ij}|)}\,(1-t)^{|S_{ij}|}\,,\]
and the similar sum $\Sigma_{ij}'$. We can now evaluate the first bracket in the right-hand side of (\ref{split2}):
\[\stacksum{(w,T)\in f^{-1}(\sigma)}{T_{11}=\ldots=T_{1n}=\emptyset}t^{\frac{1}{2}(\ell(w)-\ell(wT)-|T|)}\,(1-t)^{|T|}=\Sigma_{\lambda_1,1}\ldots\Sigma_{ij}'\ldots\Sigma_{ij}\ldots\Sigma_{1,\lambda_1'}'\,.\]
In fact, we first write the sum in the left-hand side as an iterated sum, which factors in the way shown above because $\Sigma_{ij}$ only depends on $w_{ij}[1,\lambda_i']=C_{i,j+1}[1,\lambda_i']$ (rather than the whole $w_{ij}$), by Proposition \ref{p2cols0}.

We conclude the proof by calculating the sum $\Sigma_{ij}$, the calculation for $\Sigma_{ij}'$ being similar but simpler. For simplicity, let $d:=\lambda_i'$, $w=w_{ij}$, $C:=C_{i,j+1}[1,\lambda_i']$, and $D:=C_{ij}$. Assume that $C$ differs from $D':=D(r_1,\overline{\jmath})\ldots (r_p,\overline{\jmath})$ with $1\le r_1<\ldots <r_p<j$ (possibly $p=0$) only in position $j$. Let $\Gamma_{ij}^r=\Delta\Delta'$, where
\begin{align*}
&
\begin{array}{llllll}\Delta:=(\!\!\!\!\!&(j,d+1),&(j,d+2),&\ldots,&(j,n),\\&(j,\overline{\jmath}),\\
&(j,\overline{n}),&(j,\overline{n-1}),&\ldots,&(j,\overline{d+1})\,)\,,\end{array}\\
&\;\:\Delta':=((j-1,\overline{\jmath}),\ldots,(2,\overline{\jmath}),(1,\overline{\jmath}))\,.
\end{align*}
Correspondingly, the chains $S_{ij}$ split into a head $S$, which can vary, and a fixed tail 
\[S':=((r_p,\overline{\jmath}),\ldots, (r_1,\overline{\jmath}))\,.\]
We have
\[\Sigma_{ij}=t^e(1-t)^p\stacksum{S\::\:(w,S)\in{\mathcal A}^r(\Delta)}{wS(j)=D'(j)} t^{\frac{1}{2}(\ell(wS)-\ell(w)-|S|)}\,(1-t)^{|S|}\,,\]
where $e:=\frac{1}{2}(\ell(wSS')-\ell(wS)-p)$. By Proposition \ref{p2cols0}, the sum in the right-hand side is 
\[t^{N_{C(j),D'(j)}(D[j+1,d])}(1-t)\,;\]
 note that this sum is missing when $D'=C$, which is another possibility. The exponent $e$ is calculated based on (\ref{diflengths}). 
\end{proof}


\begin{thebibliography}{10}

\bibitem{laaapm}
W.~Adamczak and C.~Lenart.
\newblock The alcove path model and {Y}oung tableaux, 2009.
\newblock {\tt math.albany.edu/math/pers/lenart}.

\bibitem{aakfpk}
E.~Ardonne and R.~Kedem.
\newblock Fusion products of {K}irillov-{R}eshetikhin modules and fermionic
  multiplicity formulas.
\newblock {\em J. Algebra}, 308:270--294, 2007.

\bibitem{asasem}
S.~Assaf.
\newblock A combinatorial proof of {LLT} and {M}acdonald positivity, 2008.
\newblock {\tt http://www-math.mit.edu/\~{}sassaf}.

\bibitem{gallsg}
S.~Gaussent and P.~Littelmann.
\newblock {LS}-galleries, the path model and {MV}-cycles.
\newblock {\em Duke Math. J.}, 127:35--88, 2005.

\bibitem{gahaha}
I.~Grojnowski and M.~Haiman.
\newblock Affine {H}ecke algebras and positivity of {LLT} and {M}acdonald
  polynomials, 2007.
\newblock {\tt http://math.berkeley.edu/\~{}mhaiman}.

\bibitem{hhlcfm}
J.~Haglund, M.~Haiman, and N.~Loehr.
\newblock A combinatorial formula for {M}acdonald polynomials.
\newblock {\em J. Amer. Math. Soc.}, 18:735--761, 2005.

\bibitem{humrgc}
J.~E. Humphreys.
\newblock {\em Reflection {G}roups and {C}oxeter {G}roups}, volume~29 of {\em
  Cambridge Studies in Advanced Mathematics}.
\newblock Cambridge University Press, Cambridge, 1990.

\bibitem{kancgr}
M.~Kashiwara and T.~Nakashima.
\newblock Crystal graphs for representations of the {$q$}-analogue of classical
  {L}ie algebras.
\newblock {\em J. Algebra}, 165:295--345, 1994.

\bibitem{katsfq}
S.~Kato.
\newblock Spherical functions and a {$q$}-analog of {K}ostant's weight
  multiplicity formula.
\newblock {\em Invent. Math.}, 66:461--468, 1982.

\bibitem{lassuc}
A.~Lascoux and M.-P. Sch\mbox{\"{u}}tzenberger.
\newblock Sur une conjecture de {H}. {O}. {F}oulkes.
\newblock {\em C. R. Acad. Sci. Paris \mbox{S\'e}r. I Math.}, 288:95--98, 1979.

\bibitem{laslqa}
C.~Lecouvey and M.~Shimozono.
\newblock Lusztig's {$q$}-analogue of weight multiplicity and one-dimensional
  sums for affine root systems.
\newblock {\em Adv. Math.}, 208:438--466, 2007.

\bibitem{lenhlp}
C.~Lenart.
\newblock Hall-{L}ittlewood polynomials, alcove walks, and fillings of {Y}oung
  diagrams, {I}.
\newblock {\tt arXiv:math.CO/0804.4715}.

\bibitem{lencfm}
C.~Lenart.
\newblock On combinatorial formulas for {M}acdonald polynomials.
\newblock {\em Adv. Math.}, 220:324--340, 2009.

\bibitem{lapawg}
C.~Lenart and A.~Postnikov.
\newblock Affine {W}eyl groups in {$K$}-theory and representation theory.
\newblock {\em Int. Math. Res. Not.}, pages 1--65, 2007.
\newblock Art. ID rnm038.

\bibitem{lapcmc}
C.~Lenart and A.~Postnikov.
\newblock A combinatorial model for crystals of {K}ac-{M}oody algebras.
\newblock {\em Trans. Amer. Math. Soc.}, 360:4349--4381, 2008.

\bibitem{litlrr}
P.~Littelmann.
\newblock {A Littlewood-Richardson rule for symmetrizable Kac-Moody algebras}.
\newblock {\em Invent. Math.}, 116:329--346, 1994.

\bibitem{litpro}
P.~Littelmann.
\newblock Paths and root operators in representation theory.
\newblock {\em Ann. of Math. {\rm (2)}}, 142:499--525, 1995.

\bibitem{litcsf}
D.~Littlewood.
\newblock On certain symmetric functions.
\newblock {\em Proc. London Math. Soc. (3)}, 11:485--498, 1961.

\bibitem{lusscq}
G.~Lusztig.
\newblock Singularities, character formulas, and a {$q$}-analog of weight
  multiplicities.
\newblock In {\em Analysis and topology on singular spaces, II, III (Luminy,
  1981)}, volume 101 of {\em Ast\'erisque}, pages 208--229. Soc. Math. France,
  Paris, 1983.

\bibitem{macsfg}
I.~Macdonald.
\newblock {\em Spherical functions on a group of {$p$}-adic type}.
\newblock Ramanujan Institute, Centre for Advanced Study in
  Mathematics,University of Madras, Madras, 1971.
\newblock Publications of the Ramanujan Institute, No. 2.

\bibitem{macsft}
I.~Macdonald.
\newblock Schur functions: theme and variations.
\newblock In {\em S\'eminaire Lotharingien de Combinatoire (Saint-Nabor,
  1992)}, volume 498 of {\em Publ. Inst. Rech. Math. Av.}, pages 5--39. Univ.
  Louis Pasteur, Strasbourg, 1992.

\bibitem{macopa}
I.~Macdonald.
\newblock Orthogonal polynomials associated with root systems.
\newblock {\em S\'em. Lothar. Combin.}, 45:Art.\ B45a, 40 pp. (electronic),
  2000/01.

\bibitem{narkfp}
K.~Nelsen and A.~Ram.
\newblock Kostka-{F}oulkes polynomials and {M}acdonald spherical functions.
\newblock In {\em Surveys in combinatorics, 2003 (Bangor)}, volume 307 of {\em
  London Math. Soc. Lecture Note Ser.}, pages 325--370. Cambridge Univ. Press,
  Cambridge, 2003.

\bibitem{ramawh}
A.~Ram.
\newblock Alcove walks, {H}ecke algebras, spherical functions, crystals and
  column strict tableaux.
\newblock {\em Pure Appl. Math. Q.}, 2:963--1013, 2006.

\bibitem{raycfm}
A.~Ram and M.~Yip.
\newblock A combinatorial formula for {M}acdonald polynomials.
\newblock {\tt arXiv:math/0803.1146}.

\bibitem{schghl}
C.~Schwer.
\newblock Galleries, {H}all-{L}ittlewood polynomials, and structure constants
  of the spherical {H}ecke algebra.
\newblock {\em Int. Math. Res. Not.}, pages Art. ID 75395, 31, 2006.

\bibitem{stekfp}
J.~Stembridge.
\newblock Kostka-{F}oulkes polynomials of general type.
\newblock {\tt www.math.lsa.umich.edu/\~{}jrs}.
\newblock Lecture notes for the {\em Generalized Kostka Polynomials Workshop},
  American Institute of Mathematics, July 2005.

\end{thebibliography}

\end{document}